\documentstyle[a4,11pt,epsf, amssymb]{article}
\begin{document}
\renewcommand{\theequation}{\arabic{section}.\arabic{equation}}
\title{\bf A Constructive Approach to the Estimation of  Dimension Reduction Directions}
\author{\\ \normalsize
Yingcun  Xia
\\
\
\\
National University of Singapore  }
\date{}
\def\beginn{\begin{eqnarray*}}
\def\endn{\end{eqnarray*}}
\def\beginy{\begin{eqnarray}}
\def\endy{\end{eqnarray}}
\def\n{\nonumber}
\newtheorem{Theorem}{Theorem}[section]
\newtheorem{Example}[Theorem]{Example}
\newtheorem{Lemma}[Theorem]{Lemma}
\newtheorem{Note}[Theorem]{Note}
\newtheorem{Proposition}[Theorem]{Proposition}
\newtheorem{Corollary}[Theorem]{Corollary}
\newtheorem{Remark}[Theorem]{Remark}
\def\ditem{\vspace{-0.1cm} \item}

\def\Cov{\mbox{\rm Cov}}
\def\Var{\mbox{\rm Var}}
\def\YX{_{Y|X}}
\def\argmin{\mbox{\rm arg}\min}
\def\b{{\mathbf b}}
\def\pperp{\perp\hspace{-.25cm}\perp}
\def\t{{\hspace{-0.05cm}\top}}
\def\Bt{B^\top}
\def\E{{\cal E}}
\def\D{{\cal D}}
\def\N{{\cal N}}
\def\O{{\cal O}}
\def\V{{{\cal V}}}
\def\R{{\Bbb R}}
\def\btd{\bigtriangledown}
\def\btdt{\bigtriangledown^\t\hspace{-0.1cm}}
\def\dfor{{\qquad \mbox{\rm for} \quad}}
\def\C{{\cal C}}
\def\rt{\raisebox{-1.5ex}[0pt]}
\def\bbar{{b\hspace{-0.1cm}\bar{} \ }}
\def\diag{\mbox{\rm diag}}
\def\tg{\tilde g}
\def\eqc{&\hspace{-0.35cm}=&\hspace{-0.35cm}}
\def\deqc{&\hspace{-0.35cm}\stackrel{def}=&\hspace{-0.35cm}}

\maketitle

\rm

\begin{abstract}

\baselineskip1.9em

{\small In this paper, we propose two new methods to estimate the
dimension-reduction directions of the central subspace (CS) by
constructing a regression model such that the directions are all
captured in the regression mean. Compared with the inverse
regression estimation methods (e.g. Li, 1991, 1992; Cook and
Weisberg, 1991), the new methods require no strong assumptions on
the design of covariates or the functional relation between
regressors and the response variable, and have better performance
than the inverse regression estimation methods for finite samples.
Compared with the direct regression estimation methods (e.g.
H\"ardle and Stoker, 1989; Hristache, Juditski, Polzehl and
Spokoiny, 2001; Xia, Tong, Li and Zhu, 2002), which can only
estimate the directions of CS in the regression mean, the new
methods can detect the directions of CS exhaustively. Consistency
of the estimators and the convergence of corresponding algorithms
are proved. }

\vspace{.5cm}


\

\noindent {\it short title}: Constructive Dimension Reduction

\

\noindent {\it AMS 200 Subject Classifications}: Primary 62G08;
Secondary 62G09, 62H05.

\end{abstract}

\newpage

\vspace{.8cm}

\baselineskip1.8em

\setcounter{equation}{0}

\section{Introduction}

Suppose $ X $ is a random vector in $ {\Bbb R}^p $ and $ Y $ is a
univariate random variable. Let $ B_0 = (\beta_{01}, \cdots,
\beta_{0q}) $ denote a  $ p\times q $ orthogonal matrix with $q
\le p$, i.e. $B_0^\t B_0 = I_q $, where $ I_q $ is a $q\times q$
identity matrix. Given $ B_0^\t X $, if $ Y $ and $ X $ are
independent, i.e. $ Y\pperp X| B_0^\t X $, then the space spanned
by the column vectors $ \beta_{01}, \beta_{02}, \cdots, \beta_{0q}
$, $ {{\cal S}}(B_0) $, is called the dimension reduction space.
If all the other dimension reduction spaces include $ {{\cal
S}}(B_0) $ as their subspace, then $ {{\cal S}}(B_0) $ is the
so-called central dimension reduction subspace (CS); see Cook
(1998). The column vectors $ \beta_{01}, \beta_{02}, \cdots,
\beta_{0q} $ are called the CS directions. Dimension reduction is
a fundamental statistical problem both in theory and in practice.
See Li (1991, 1992) and Cook (1998) for more discussion. If the
conditional density function of $ Y $ given $ X$ exists, then the
definition is equivalent to the conditional density function of $
Y | X$ being the same as that of $ Y|\Bt_0 X $ for all possible
values of $ X $ and $ Y$, i.e.
\beginy
f_{_{Y|X}}(y|x) = f_{_{Y|B_0^\t X}}(y|B_0^\t x).
\label{definition}
\endy
Other alternative definitions for the dimension reduction space
can be found in the literature.

In the last decade or so, a series of papers (e.g. H\"ardle and
Stocker, 1989; Li, 1991; Cook and Weisberg, 1991; Samarov, 1993;
Hristache, Juditski, Polzehl and Spokoiny, 2001; Yin and Cook,
2002; Xia, Tong, Li and Zhu, 2002; Cook and Li, 2002; Li, Zha and
Chiaromonte, 2004; Lue, 2004) have considered issues related to
the dimension reduction problem, with the aim of estimating the
dimension reduction space and relevant functions. The estimation
methods in the literature can be classified into two groups:
inverse regression estimation methods (e.g. SIR, Li, 1991 and
SAVE, Cook and Weisberg, 1991) and direct regression estimation
methods (e.g. ADE, H\"ardle and Stoker, 1991 and MAVE of Xia,
Tong, Li and Zhu 2002). The inverse regression estimation methods
are computationally easy and are widely used as an initial step in
data mining, especially for large data sets. However, these
methods have poor performance in finite samples and need strong
assumptions on the design of covariates. The direct regression
estimation methods have much better performance for finite samples
than the inverse regression estimations. They need no strong
requirements on the design of covariates or on the response
variable. However, the direct regression estimation methods cannot
find the directions in CS exhaustively, such as those in the
conditional variance.

None of the methods mentioned above use the definitions directly
in searching for the central space. As a consequence, they fail in
one way or another to estimate CS efficiently. A straightforward
approach in using definition (\ref{definition}) is to look for $
B_0 $ in order to minimize the difference between those two
conditional density functions. The conditional density functions
can be estimated using nonparametric smoothers. Obviously, this
approach is not efficient in theory due to the ``curse of
dimensionality'' in nonparametric smoothing. In calculations, the
minimization problem is difficult to implement. People have
observed that the CS in the regression mean function, i.e. the
central mean space (CMS), can be estimated much more efficiently
than the general CS. See, for example, Yin and Cook (2002), Cook
and Li (2002) and Xia, Tong, Li and Zhu (2002). Motivated by this
observation, one can construct a regression model such that the CS
coincides with the CMS space in order to reduce the difficulty of
estimation. In this paper, we first construct a regression model
in which the conditional density function $ f_{Y|X}(y|x) $ is
asymptotically equal to the conditional mean function. Then, we
apply the methods of searching for the CMS to the constructed
model. Based on the discussion above, this constructive approach
is expected to be more efficient than the inverse regression
estimation methods for finite samples, and can detect the CS
directions exhaustively.

In the estimation of dimension reduction space,  most methods need
in one way or another to deal with nonparametric estimation. In
terms of nonparametric estimation, the inverse regression
estimation methods employ a nonparametric regression of $ X$ on  $
Y $ while the direct regression estimation methods employ a
nonparametric regression of $ Y $ on $ X $.  In contrast to
existing methods, the methods in this paper search for CS from
both sides by investigating conditional density functions. A
similar idea appeared in Yin and Cook (2005) for a general
single-index model. To overcome the difficulties of calculation,
we propose two algorithms in this paper using a similar idea to
Xia, Tong, Li and Zhu (2002). The algorithm solves the
minimization problem in the method by treating it as two separate
quadratic programming problems, which have simple analytic
solutions and can be calculated quite efficiently. The convergence
of the algorithm can be proved. Our constructive approach can
overcome the disadvantages both in inverse regression estimations,
requiring a symmetric design for explanatory variables, and also
the disadvantage in direct regression estimation, of not finding
the CS directions exhaustively. Simulations  suggest that the
proposed methods have very good performance for finite samples and
are able to estimate the CS directions  in very complicated
structures. Applying the proposed methods to two real data sets,
some useful patterns have been observed, based on the estimations.

To estimate the CS, we need to estimate the directions $ B_0 $ as
well as the  dimension $ q $ of the space. In this paper, however,
we focus on the estimation of the directions by assuming that $ q
$ is known.

\setcounter{equation}{0}

\section{Estimation methods}

As discussed above, the direct regression estimations have good
performance for finite samples. However, it cannot detect
exhaustively the CS directions   in complicated structures.
Motivated by these  facts, our strategy is to construct a
semiparametric regression model such that all the CS directions
are captured in the regression mean function. As we can see from
(\ref{definition}), all the directions can  be captured in the
conditional density function. Thus, we will construct a regression
model such that the conditional density function is asymptotically
equal to the regression mean function.

The primary step is thus to construct an estimate for the
conditional density function. Here, we use the idea of the
``double-kernel'' local linear smoothing method studied in Fan et
al (1996) for the estimation. Consider $ H_b(Y-y) $ with $ y $
running through all possible values, where $ H(v) $ is a symmetric
density function, $ b > 0 $ is a bandwidth and $ H_b(v) = b^{-1}
H(v/b)$. If $ b \to 0 $ as $ n \to \infty $, then from
(\ref{definition}) we have
\beginn
m_b(x, y) \stackrel{def}=  E(H_b(Y-y)|X=x) = E(H_b(Y-y)|B_0^\t X=
B_0^\t x) \to f_{Y|B_0^\t X}(y|B_0^\t x).
\endn
See Fan et al (1996). The above equation indicates  that all the
directions can be captured by the conditional mean function $
m_b(x, y) $ of $ H_b(Y-y) $ on $ X = x $ with $ x$ and $ y $
running through all possible values. Now, consider a regression
model nominally for $H_b(Y-y) $ as
\beginn
  H_b(Y-y) = m_b(X, y) + \varepsilon_{_b}(y|X),
\endn
where $ \varepsilon_{_b}(y|X) = H_b(Y-y) - E(H_b(Y-y)| X)  $ with
$ E \varepsilon_{_b}(y|X) = 0 $. Let $ g_b(B_0^\top x, y) =
E(H_b(Y-y)|B_0^\t X= B_0^\t x) $. If (\ref{definition}) holds,
then $ m_b(x, y) = g_b(B_0^\top x, y) $. The model can be written
as
\beginy
  H_b(Y-y) = g_b(B_0^\top X, y) + \varepsilon_{_b}(y|X).  \label{model}
\endy
As $ b \to 0 $, we have $ g_b(B_0^\top x, y) \to f_{Y|B_0^\top
X}(y|B_0^\top x) $. Thus, the directions $ B_0 $ defined in
(\ref{definition}) are all captured in the regression mean
function in model (\ref{model}) if $ y $ runs through all possible
values.

Based on model (\ref{model}), we propose two  methods to estimate
the directions.  One of the methods is a combination of the outer
product of gradients (OPG) estimation method (H\"ardle, 1991;
Samarov, 1993; Xia, Tong, Li and Zhu, 2002) with the
``double-kernel'' local linear smoothing method (Fan et al, 1996).
The other one is a combination of the minimum average
(conditional) variance estimation (MAVE) method (Xia, Tong, Li and
Zhu, 2002) with the ``double-kernel'' local linear smoothing
method. The structure adaptive weights in Hristache, Juditski and
Spokoiny (2001) and Hristache, Juditski, Polzehl and Spokoiny
(2001) are used  in the estimations.

\subsection{Estimation based on outer products of gradients}

Consider the gradient of the conditional mean function $ m_b(x, y)
$ with respect to $ x$. If (\ref{definition}) holds, then it
follows
\beginy
  \frac{\partial m_b(x, y)}{\partial x} = \frac{\partial g_b(B_0^\top x, y)}{\partial x} = B_0 \btd
 g_b(B_0^\top x, y),  \label{grad}
\endy
where $ \btd g_{_b}(v_1, \cdots, v_q,y) = ( \btd_1 g_{_b}(v_1,
\cdots, v_q, y), \cdots, \btd_q g_{_b}(v_1, \cdots, v_q, y))^\t $
with
$$
\btd_k g_{_b}(v_1, \cdots, v_q,y) = \frac{\partial}{\partial v_k}
g_{_b}( v_1, \cdots, v_q, y), \quad k = 1, 2, \cdots, q.
$$
Thus, the directions $ B_0 $ are contained in the gradients of the
regression mean function in model (\ref{model}). One way to
estimate $ B_0 $ is by considering the expectation of the outer
product of the gradients
\beginn
  E\Big\{\Big(\frac{\partial m_b(X, Y)}{\partial x}\Big)
  \Big(\frac{\partial m_b(X, Y)}{\partial x}\Big)^\top\Big\} =
  B_0 E\{ \btd g_b(B_0^\top X, Y) \btdt g_b(B_0^\top X,
  Y)\} B^\top_0.
\endn
It is easy to see that  $ B_0$ is in the space spanned by the
first $q$ eigenvectors of the expectation of the outer products.

Suppose that $ \{(X_i, Y_i), i = 1, 2, \cdots n\} $ is a random
sample from $ (X, Y) $. To estimate the gradient $ \partial m_b(
x, y)/\partial x $, we can use the nonparametric kernel smoothing
methods.  For simplicity, we adopt the following notation scheme.
Let $ K_0(v^2) $ be a univariate symmetric density function and
define $ K(v_1, \cdots, v_d) = K_0( v_1^2 + \cdots + v_d^2) $ for
any integer $d $ and $ K_h(u) = h^{-d} K(u/h)$, where $ d $ is the
dimension of $ u $ and $ h > 0 $ is a bandwidth. Let $H_{b,i}(y) =
H_b(Y_i-y) $, where $ H(.) $ and $ b $ are defined above. For any
$ (x, y)$, the principle of the local linear smoother suggests
minimizing
\beginy
 n^{-1} \sum_{i=1}^n \Big\{ H_{b,i}(y) - a - b^\t
 (X_i-x)\Big\}^2 K_h( X_{ix} )  \label{estg}
\endy
with respect to $ a $ and $ b $ to estimate $ m_b(x, y)$ and $
\partial m_b( x, y)/\partial x $ respectively, where $
X_{ix} = X_i - x$. See Fan and Gijbels (1996) for more details.
For each pair of $ (X_j, Y_k) $, we consider the following
minimization problem
\beginy
(\hat a_{jk}, \hat b_{jk}) = \mbox{arg}\min_{a_{jk},
b_{jk}}\sum_{i=1}^n \Big[ H_{b,i} (Y_k) - a_{jk} - b_{jk}^\t
X_{ij}
    \Big]^2 w_{ij},        \label{fjkrg9}
\endy
where $ X_{ij} = X_i - X_j $ and $ w_{ij} = K_h( X_{ij} )$.  We
consider an average of their outer products
\beginn
\hat \Sigma = n^{-2}\sum_{k=1}^n  \sum_{j=1}^n \hat \rho_{jk} \hat
b_{jk} \hat b_{jk}^\t ,
\endn
where $ \hat \rho_{jk} $ is a trimming function introduced for
technical purpose to handle the notorious boundary points.  In
this paper, we adopt the following trimming scheme. For any given
point $(x, y) $, we use all observations to estimate its function
value and its gradient as in (\ref{estg}). We then consider the
estimates in a compact region of $ (x, y) $. Moreover, for those
points with too few observations around, their estimates might be
unreliable. They should not be used in the estimation of the CS
directions  and should be trimmed off. Let $\rho (\cdot )$ be any
bounded function with bounded second order derivatives on
${\mathbb{R}}$ such that $\rho (v) > 0$ if $v>\omega_0$; $\rho
(v)=0$ if $v \le \omega_{0}$ for some small $ \omega_{0}>0$. We
take $ \hat \rho_{jk} = \rho (\hat f(X_j))\rho( \hat f_{_Y}(Y_k))
$, where $ \hat f(x) $ and $ \hat f_{_Y}(y) $ are estimators of
the density functions of $ X $ and $ Y $ respectively. The CS
directions can be estimated by the first $q$ eigenvectors of $
\hat \Sigma $.

To allow the estimation to be adaptive to the structure of the
dependency of $ Y$ on $ X $, we may follow the idea of Hristache
et al (2001) and replace $ w_{ij} $ in (\ref{fjkrg9}) by
\beginn
w_{ij} = K_h( \hat \Sigma^{1/2}X_{ij}),
\endn
where $ \hat \Sigma^{1/2} $ is a symmetric matrix with $ (\hat
\Sigma^{1/2})^2 = \hat \Sigma$. Repeat the above procedure until
convergence. We call this procedure the method of outer product of
gradient based on the conditional density functions (dOPG). To
implement the estimation procedure, we suggest the following dOPG
algorithm.

\begin{enumerate}
\item[Step] 0: Set $ \hat \Sigma_{(0)} = I_p $ and $ t = 0 $.

\def\mo{\footnotesize -1}
\item[Step] 1: With $ w_{ij} = K_h( \hat \Sigma_{(t)}^{1/2}X_{ij})
$, calculate the solution to (\ref{fjkrg9})
\beginn
{a^{(t)}_{jk}\choose b^{(t)}_{jk}} &=& \Big\{\sum_{i=1}^n
   K_{h_{t}}(\hat \Sigma^{1/2}_{(t)} X_{ij})
   {1 \choose   X_{ij} } {1 \choose   X_{ij}
   }^{\t}\Big\}^{-1}  \\
   && \hspace{3cm} \times   \sum_{i=1}^n K_{h_{t}}(\hat\Sigma^{1/2}_{(t)} X_{ij})
   {1 \choose  X_{ij} } H_{b_{t},i}(Y_k),
\endn
where  $ h_{t} $ and $ b_{t} $ are bandwidths (details are given
in (\ref{band0}) and (\ref{band1}) below).

\def\mo{\footnotesize -1}
\item[Step] 2: Define $ \rho^{(t)}_{jk} = \rho(\tilde
f^{(t)}(X_j))\rho(  \tilde f^{(t)}_{_Y}(Y_k))   $ with
$$
\tilde f^{(t)}_{_Y}(y) = n^{-1}\sum_{i=1}^n H_{b_t,i}(y),\quad
\tilde f^{(t)}(x) = (n \tilde \mu )^{-1} h_t^p
\prod_{\lambda^{(t)}_k > h_t} \frac{\lambda^{(t)}_k}{h_t}
\sum_{i=1}^n K_{h_t}( \hat \Sigma^{1/2}_{(t)} X_{ix}),
$$
where $ \lambda^{(t)}_k, k =1, \cdots, p, $ are the eigenvalues of
$ \hat \Sigma^{1/2}_{(t)} $ and $ \tilde \mu = \int K_0(
\sum\limits_{\lambda^{(t)}_k
> h_t} \hspace{-0.2cm} v_k^2)$ $\prod_{\lambda^{(t)}_k> h_t} dv_k$. Calculate the average of outer products
$$
\hat \Sigma_{(t+1)} = n^{-2} \sum_{j,k=1}^n \rho^{(t)}_{jk}
b^{(t)}_{jk} (b^{(t)}_{jk})^{\t}.
$$

\item[Step] 3:  Set $ t := t+1$. Repeat Steps 1 and 2 until
convergence. Denote the final value of $ \hat \Sigma_{(t)} $ by $
 \Sigma_{(\infty)} $. Suppose the eigenvalue decomposition of $  \Sigma_{(\infty)} $ is
$ \Gamma diag(\lambda_1, \cdots, \lambda_p) \Gamma^\top $, where $
\lambda_1 \ge  \cdots \ge \lambda_p $. Then the estimated
directions are the first $ q $ columns of $ \Gamma $, denoted by $
\hat B_{dOPG}$.

\end{enumerate}

In the dOPG algorithm, $ \tilde f_Y^{(t)}(y) $ and $ \tilde
f^{(t)}(x) $, $ t > 0 $, are the estimators of the density
functions of $Y $ and $ B_0^\top X $ respectively. A justification
is given in the proof of Theorem \ref{mainOPG} in Section 6.2. In
calculations, the usual stopping criterion can be used. For
example, if the largest singular value of $ \hat \Sigma_{(t)} -
\hat \Sigma_{(t+1)} $ is smaller than $ 10^{-6}$ then we stop the
iteration and take $ \hat \Sigma_{(t+1)} $ as the final estimator.
The eigenvalues of $ \Sigma_{(\infty)} $ can be used to determine
the dimension of the CS. However, we will not go into the details
on this issue in this paper. In practice, we may need to
standardize $ X_i = (X_{i1}, \cdots, X_{ip})^\t $ by setting $ X_i
:= S_{_X}^{-1/2}(X_i - \bar X)$ and standardize $ Y_i $ by setting
$ Y_i := (Y_i-\bar Y)/\sqrt{s_{_Y}} $, where $ \bar X = n^{-1}
\sum_{i=1}^n X_i $ and $ S_{_X} = n^{-1}\sum_{i=1}^n ( X_{i} -
\bar X) ( X_{i} - \bar X)^\t $, $ \bar Y = n^{-1} \sum_{i=1}^n Y_i
$ and $  s_{_Y} = n^{-1}\sum_{i=1}^n ( Y_i - \bar Y )^2 $. Then
the estimated CS directions  are the first $ q $ columns of  $
\Gamma S_{_X}^{-1/2}$.

\subsection{MAVE based on conditional density function}

Note that  if (\ref{definition}) holds, then  the gradients $
\partial m_b(x, y)/\partial x $ at all  $ (x, y) $ are in
a common $q$-dimensional subspace as shown in equation
(\ref{grad}). To use  this observation, we can replace $ b $ in
(\ref{estg}), which is an estimate of the gradient, by $ B d(x,y)
$ and have the following local linear approximation
\beginn
  n^{-1} \sum_{i=1}^n \{ H_{b,i}(y) - a - d^\top B^\top (X_i - x)
  \}^2 K_h(X_{ix}),
\endn
where $ d = d(x, y) $ is introduced to take the role of $ \btd
g_b(B_0^\top x, y) $ in (\ref{grad}). Note that the above weighted
mean of squares is the local approximation errors of $ H_{b,i}(y)
$ by a hyperplane with the normal vectors in a common space
spanned by $ B $. Since $ B $ is common for all $ x $ and $ y $,
it should be estimated with aims to minimize the approximation
errors for all possible $ X_j $ and $ Y_k $. As a consequence, we
propose to estimate $B_0$ by minimizing
\beginy
 n^{-3} \sum_{k=1}^n \sum_{j=1}^n \hat \rho_{jk} \sum_{i = 1}^n \{H_{b,i}(Y_k)
- a_{jk} - d^{\t}_{jk}\Bt X_{ij} \}^2 w_{ij} \label{mini2}
\endy
with respect to $ a_{jk}, d_{jk}=(d_{jk1}, \cdots, d_{jkq})^{\t},
j,k=1,...,n $ and $B: \Bt B = I_q $, where $ \hat \rho_{jk} $ is
defined above. This estimation procedure is similar to the minimum
average (conditional) variance estimation method (Xia, Tong, Li
and Zhu, 2002). Because the method is based on the conditional
density functions, we call it the minimum average (conditional)
variance estimation based on the conditional density functions
(dMAVE).

The minimization problem in (\ref{mini2}) can be solved by fixing
$ (a_{jk}, d_{jk}), j,k=1,...,n, $ and $B $ alternatively. As a
consequence, we need to solve two quadratic programming problems
which have simple analytic solutions. For any matrix $ B =
(\beta_1, \cdots, \beta_d) $, we define operators $ \ell(.) $ and
$ {{\cal M}}(.) $ respectively as
\beginn
\ell(B) = (\beta_1^\t, \cdots, \beta_d^\t)^\t \qquad
\mbox{and}\qquad {{\cal M}}(\ell(B)) = B.
\endn
We propose the following dMAVE algorithm to implement the
estimation.

{
\begin{enumerate}
\ditem[Step] 0: Let $ B_{(1)} $ be an initial estimator of the CS
directions. Set $ t = 1 $.

\ditem[Step] 1: Let $ B = B_{(t)}$, calculate the solutions of $
(a_{jk}, d_{jk}), j,k=1,...,n, $ to the minimization problem in
(\ref{mini2})
\beginn
{a^{(t)}_{jk}\choose d^{(t)}_{jk}} &=& \Big\{\sum_{i=1}^n
   K_{h_{t}}(B_{(t)}^{\t}X_{ij})
   {1 \choose  B_{(t)}^{\t} X_{ij} } {1 \choose  B_{(t)} ^{\t} X_{ij}
   }^{\t}\Big\}^{-1} \\
 && \hspace{4cm} \times   \sum_{i=1}^n K_{h_{t}}(B_{(t)}^{\t}X_{ij})
   {1 \choose  B_{(t)}^{\t}X_{ij} } H_{b_{t},i}(Y_k),
\endn
where  $ h_{t} $ and $ b_{t} $ are two bandwidths (details are
discussed below).

\ditem[Step] 2: Let $ \rho_{jk}^{(t)} =\rho(\hat
f_{B_{(t)}}({X_j}))\rho (\hat f^{(t)}_{_Y}(Y_k) ) \ $ with  $ \
\hat f^{(t)}_{_Y}(y) = $  $ n^{-1}\sum_{i=1}^n H_{b_{t},i}(y) $
and  $ \hat f_{B_{(t)}}(x) = n^{-1} \sum_{i=1}^n
K_{h_{t}}(B^\top_{(t)} X_{ix}) $.  Fixing $ a_{jk} = a_{jk}^{(t)}
$ and $ d_{jk} = d^{(t)}_{jk} $, calculate the solution of $ B $
or $ \ell(B) $ to (\ref{mini2})
\beginn
\b^{(t+1)} &=& \Big\{\sum_{k,j,i=1}^n \rho_{jk}^{(t)}
            K_{h_{t}}( B_{(t)}^{\t} X_{ij})
        X^{(t)}_{ijk} (X^{(t)}_{ijk})^{\t}
       \Big\}^{-1} \\
 &&  \hspace{3cm} \times \sum_{k,j,i=1}^n \rho_{jk}^{(t)} K_{h_{t}}(
B_{(t)}^{\t} X_{ij})
        X^{(t)}_{ijk}\{H_{b_{t},i}(Y_k) - a^{(t)}_{jk}\},
\endn
where  $ X_{ijk}^{(t)} = d^{(t)}_{jk} \otimes X_{ij} $.

\item[Step] 3: Calculate $ \Lambda_{(t+1)} = \{ {{\cal
M}}(\b^{(t+1)})\}^\top {{\cal M}}(\b^{(t+1)}) $ and $ B_{(t+1)} =
{{\cal M}}(\b^{(t+1)}) \Lambda_{(t+1)}^{-1/2}$. Set $ t:= t +1 $
and go to Step 1.

\item[Step] 4: Repeat steps 1--3 until convergence. Let $
B_{(\infty)} $ be the final value of $ B_{(t)} $. Then our
estimators of the directions are the columns in $ B_{(\infty)} $,
denoted by $ \hat B_{dMAVE} $.

\end{enumerate}

}

The dMAVE algorithm needs a consistent initial estimator in Step 0
to guarantee its theoretical justification. In the following, we
use the first iteration estimator of dOPG, the first $ q$
eigenvector of $ \hat \Sigma_{(1)} $, as the initial value.
Actually, any initial estimator that satisfies (\ref{rate0}) can
be used and Theorem \ref{mainMAVE} will hold. Similar to dOPG, the
standardization procedure can be carried out for dMAVE in
practice. The stopping criterion for dOPG can also be used here.

Note that the estimation in the procedure is related with
nonparametric estimations of conditional density functions.
Several bandwidth selection methods are available for the
estimation. See, e.g. Silverman (1986), Scott (1992) and Fan et al
(1996). Our theoretical verification of the convergence for the
algorithms requires some constraints on the bandwidths although we
believe these constraints can be removed with more complicated
technical proofs. To ensure the requirements on bandwidths can be
satisfied, after standardizing the variables we use the following
bandwidths in our calculations. In the first iteration, we use
slightly larger bandwidths than the optimal ones in terms of MISE
as
\beginy
h_0 = c_0 n^{-\frac1{p_0+6}}, \qquad b_0 = c_0 n^{-\frac1{p_0+5}},
\label{band0}
\endy
where $ p_0 = \max(p, 3) $. Then we reduce the bandwidths in each
iteration as
\beginy
h_{t+1} = \max\{r_n h_t,  c_0 n^{-\frac1{q+4}}\}, \qquad b_{t +1}
= \max\{r_n b_t, c_0n^{-\frac1{q+3}}, c_0 n^{-\frac15}\}
\label{band1}
\endy
for $ t \ge 0 $, where $ r_n = n^{-1/(2(p_0+6))} $, $ c_0 = 2.34$
as suggested by Silverman (1986) if the Epanechnikov kernel is
used. Here, the bandwidth $ b $ is selected  smaller than $ h$
based on simulation comparisons.

Fan and Yao (2003, p.337) proposed a method, called the profile
least-squares estimation, for the single-index model and its
variants by solving a similar minimization problem as in
(\ref{mini2}). The method is also possible to be used here for the
estimation of $ B_0 $ in (\ref{model}).

\setcounter{equation}{0}

\section{Asymptotic results}

To exclude the trivial cases, we assume that $ p > 1 $ and $ q \ge
1 $. Let $ f_0(y|v_1, \cdots, v_q) $, $ f_0(v_1,\ \cdots,\ v_q) $
and $ f_{_Y}(y) $ be the (conditional) density functions of $ Y|
B_0^\top X $, $ B_0^\top X$ and $ Y $ respectively. Let \ $
\rho_0(x, y) = \rho(f_0(B_0^\top x)) \rho(f_{_Y}(y)) $, \ $
 \btd f_0(y| v_1, \cdots, v_q) \ = $ $ \ (\partial f_0(y|v_1, \cdots, v_q
 )/\partial v_1, \cdots, \partial f_0(y|v_1, \cdots, v_q
 )/\partial v_q)^\top
$, $ \mu_B(u) = E(X|B^\top X = u) $ and $ w_B(u) = E\{XX^\top|B^\t
X=u\}$. For any matrix $ A$, let $ |A| $ denote its largest
singular value, which is same as the Eculidean norm if $ A $ is a
vector. Let $ \tilde B_0: p\times (p-q) $ be such that $ (B_0,
\tilde B_0)^\top (B_0, \tilde B_0) = I_p $. We need the following
conditions for (\ref{definition}) to prove our theoretical
results.

\begin{enumerate}
\item[(C1)]  [Design of $ X$] The density function $f(x) $ of $ X
$ has bounded  second order derivatives on $ \R^p $; $ E|X|^r <
\infty $ for some $ r > 8 $; functions $ \mu_B(u) $ and $ w_B(u) $
have bounded derivatives with respect to $ u $ and $ B $  for  $ B
$ in a small neighbor of $ B_0 $: $ |B-B_0| \le \delta $ for some
$ \delta> 0$.

\vspace{-0.2cm}

\item[(C2)]  [Conditional density function] The conditional
density functions $ f_{Y| X}(y|x) $ and  $ f_{Y|B^\top X}$ $(y|u)
$ have bounded fourth order derivatives with respect to $ x$, $ u
$ and $ B$  for $ B $ in a small neighbor of $ B_0 $; the
conditional density function of $ f_{\tilde B_0^\top X, Y|
B_0^\top X}(u, y|v) $ and  $ \int |\btd f_0(y|u)| dy $ are bounded
for all $ u, y $ and $ v$.

\vspace{-0.2cm}

\item[(C3)]  [Efficient dimension] Matrix $ M_0 = \int \rho_0(x,
y) \btd f_0(y|B_0^\t x) \btdt f_0(y| B_0^\t x)  f(x) $ $ f_{_Y}(y)
dxdy$ has full rank $ q$.

\vspace{-0.2cm}

\item[(C4)]  [Kernel functions] $K_0(v^2)$ and $ H(v) $ are two
symmetric univariate density functions with bounded second order
derivatives and compact supports.

\vspace{-0.2cm}

\item[(C5)]  [Bandwidths for consistency] Bandwidths $ h_0 = c_1
n^{-r_h} $ and $ b_0 = c_2 n^{-r_b} $ where $ 0 < r_h, r_b \le
1/(p_0+6)$, $ p_0 = \max\{p, 3\} $. For $ t \ge 1 $, $ h_t =
\max\{ r_n h_{t-1}, \hbar  \}$ and $ b_t  = \max\{ r_n b_{t-1},
\bbar \} $ where $ r_n = n^{-r_h/2} $, $ \hbar = c_3 n^{-r_h'},
\bbar = c_4 n^{-r_b'} $ with $ 0 < r_h', r_b' \le 1/(q+3) $, and $
c_1, c_2, c_3, c_4 $ are constants.


 \vspace{-0.2cm}

\end{enumerate}

In (C1), the finite moment requirement for $|X|$ can be removed if
we adopt the trimming scheme of H\"ardle et al (1993). However, as
noticed in Delecroix et al (2004), this scheme caused some
technical problems in the proofs. Based on assumptions (C2) and
(C4), the smoothness of $ g_b(u, y) $ is implied. Lower order of
smoothness is sufficient if we are only interested in the
 estimation consistency. The second order differentiable
requirement in (C4) can ensure the Fourier transformations of the
kernel functions being absolutely integrable; see Chung (p.166,
1968). The popular kernel functions such as Epanechnikov kernel
and quadratic kernel are included in (C4). The Gaussian kernel can
be used with some modifications to the proofs. Condition (C3)
indicates that the dimension $ q $ cannot be further reduced. For
ease of exposition, we further assume that $\mu _{0H} =\int H(v)
dv = 1, \mu _{2H} =\int v^2 H(v) dv = 1$, $\mu_{0q} = \int K(v_1,
\cdots,  v_q)  dv_1\cdots dv_q =1 $ and $ \mu_{2q} = \int K(v_1,
\cdots, v_q) v_1^2 dv_1\cdots dv_q =1$; otherwise, we take
$H(v):=H(v/\tau_{2H}^{1/2})/\tau _{2H}^{1/2}$ and  $ K(v_1,
\cdots, v_q) = \mu_{0q}^{-1} K(v_1/\sqrt{\mu_{2q}}, \cdots, $
$v_q/\sqrt{\mu_{2q}}) $ $/\sqrt{\mu_{2q}} $.  The bandwidths
satisfying  (C5) can be found easily. For example, the bandwidths
given in (\ref{band0}) and (\ref{band1}) satisfy the requirements.
Actually, a wider range of bandwidths can be used; see the proofs.
Let $ \nu_B(x) = \mu_B(B^\top x) - x$,  $ \bar w_B(x) = w_B(B^\top
x) - \mu_B( B^\top x)\mu_B^\top( B^\top x) $ and $ f_0(x) =
f_0(B_0^\top x) $. For any square matrix $ A$, $ A^{-1} $ and $
A^+$ denote the inverse (if it exists) and the Moore-Penrose
inverse matrices respectively.

\begin{Theorem}\label{mainOPG} \rm
Suppose  conditions (C1)-(C5) hold. Then we have
\beginn
 |  \hat B_{dOPG} \hat B^\top_{dOPG} - B_0 B_0^\top | = O(\hbar ^4 + \delta^2_{q\hbar  \bbar }+ \delta_{q\hbar  \bbar }
 \bbar ^4 + \delta_n^2/\bbar ^2 + n^{-1/2}
  )
\endn
in probability as $ n \to \infty $, where $ \delta_{q \hbar  \bbar
} = (n \hbar ^q \bbar /\log n)^{-1/2} $ and $ \delta_n = (\log
n/n)^{1/2}$. If $\hbar ^4 + \delta^2_{q\hbar  \bbar }+
\delta_{q\hbar  \bbar }
 \bbar ^4 + \delta_n^2/\bbar ^2 = o( n^{-1/2}) $ can be satisfied, then
\beginn
\sqrt{n} \{\ell(\hat B_{dOPG} \hat B^\top_{dOPG}B_0 ) - \ell( B_0
)\} \stackrel{D}\to N(0, W_0),
\endn
where
\beginn
 W_0 &=& Var [ \rho_0(X, Y) M_0^{-1} (\btd f_0(Y|B_0^\t X)
f_{_Y}(Y)
 - E\{\btd f_0(Y|B_0^\t X) f_{_Y}(Y) | X\} )\\
 && \hspace{8.5cm}
  \otimes ( \bar w_{_{B_0}}^{+}(X)
\nu_{_{B_0}}(X))] .
\endn

\end{Theorem}

The first part of Theorem \ref{mainOPG} indicates that $ \hat
B_{dOPG} $ is a consistent estimator of an orthogonal basis, $ B_0
Q $ with $ Q = B_0^\top \hat B_{dOPG}$,  in CS  and $ |\hat
B_{dOPG} - B_0 Q | = O(\hbar ^4 + \delta^2_{q\hbar \bbar }+
\delta_{q\hbar \bbar }
 \bbar ^4 + \delta_n^2/\bbar ^2 + n^{-1/2}
)$ in probability. See Bai et al (1991) and Xia, Tong, Li and Zhu
(2002) for alternative presentations of the asymptotic results. If
the bandwidths in (\ref{band1}) are used, then the consistency
rate is $ O(n^{-4/(q+4)+1/(q+3)}\log n + n^{-1/2}) $  in
probability. Faster consistency rate can be obtained by adjusting
the bandwidths. The convergence of the corresponding algorithm is
also implied in the proof in section 6. If $ q \le 3 $, then the
condition for the normality can be satisfied by taking
\beginn
1> r_h' > \frac18, \qquad \frac27 r_h' < r_b' < \frac12 - q r_h'.
\endn
If we use higher order polynomial smoothing, it is possible to
show that the root-$n $ consistency can be achieved for any
dimension $ q $. See, e.g.  H\"{a}rdle and Stoker (1989) and
Samarov (1993), where the higher order kernel, a counterpart of
the higher order polynomial smoother,  was used. However, using
higher order polynomial smoothers increases the difficulty of
calculations while the improvement of finite sample performance is
not substantial.

\begin{Theorem}\label{mainMAVE} \rm
If conditions (C1)-(C5) holds, then
\beginn
 |  \hat B_{dMAVE} \hat B_{dMAVE}^\top   - B_0  B_0^\top | = O\{\hbar ^4 + \delta^2_{q\hbar  \bbar }+ \delta_{q\hbar  \bbar }
 \bbar ^4 + \delta_n^2/\bbar ^2 + n^{-1/2}
  \}
\endn
 in probability as $ n \to \infty $. If $\hbar ^4 + \delta^2_{q\hbar  \bbar }+ \delta_{q\hbar  \bbar }
 \bbar ^4 + \delta_n^2/\bbar ^2 = o( n^{-1/2}) $ can be satisfied, then
\beginn
\sqrt{n} \{\ell(\hat B_{dMAVE} \hat B_{dMAVE}^\top B_0  ) - \ell(
B_0 )\} \stackrel{D}\to N(0, D_0^+ \Sigma_0 D_0^+),
\endn
where $ D_0 = \int \rho_0(x,y)  \btd f_0(y|B_0^\t x) \btdt f_0(y|
B_0^\t x)  \otimes \{\nu_{_{B_0}}(x) \nu_{_{B_0}}^\t (x) \} f_0(x)
f_{_Y}(y) dxdy $ and
\beginn
 \Sigma_0 = Var[ \rho_0(X, Y)
( \btd f_0(Y|B_0^\t X) f_{_Y}( Y) - E\{\btd f_0(Y|B_0^\t X)
f_{_Y}( Y)| X\} )\otimes \nu_{_{B_0}}(  X)].
\endn
\end{Theorem}

The proof of Theorem \ref{mainMAVE} is given in section 6. The
convergence of the dMAVE algorithm is implied in the proof.
Similar remarks on dOPG are applicable to dMAVE. Moreover, $ \hat
B_{dMAVE} $ converges to $ B_0 \tilde Q $, where $ \tilde Q $ is
determined by the initial consistent estimator of the directions.
For example, $ \tilde Q = \hat B_{(1)}^\top B_0 $ if $ B_{(1)} $
is used as the initial estimator. Similarly, the root-$n$
consistency holds for $ q\le 3 $. It is possible that the root-$n
$ consistency holds for $ q > 3 $ if higher order local polynomial
smoothing method is used. In spit of the equivalence in terms of
consistency rate for both dOPG and dMAVE, our simulations suggest
that dMAVE has better performance than dOPG in finite samples.
Theoretical comparison of efficiencies between the two methods is
not clear. In a very special case when $ q = 1 $ and the CS is in
the regression mean, Xia (2006a) proved that dMAVE is more
efficient than dOPG.

We here give some discussions about the requirements on the
distributions of $ X $ and $ Y$. If $ Y $ is discrete, we can
consider the conditional cumulative distribution functions and
have $ F_{Y|X}(y|x) = F_{Y|B_0^\top X}(y|B_0^\top x) $ when $ Y
\pperp X | B_0^\top X$ holds. Similar to (\ref{model}), we can
consider a regression model
\beginn
 I(Y < y) = G(B_0^\top X, y) + e(y|X),
\endn
where $ G(B_0^\top x, y) = E\{I(Y < y)|X=x\} = E\{ I(Y < y)|
B_0^\top X = B_0^\top x\}$ and $ e(y|X) = I(Y < y) - G(B_0^\top X,
y) $. Similar theoretical consistency results are possible to be
obtained following the same techniques developed here. If some
covariates in $ X $ are discrete, our algorithms in searching for
a consistent initial estimator will fail. However, if a consistent
initial estimator can be found by for example the methods in
Horowitz and H\"ardle  (1996) and Hristache, Juditski, Polzehl and
Spokoiny (2001) and that $ B^\top X $ has a continuous density
function for all $ B $ in a neighbor around $ B_0$, then our
theoretical results in the above theorems still hold.

\setcounter{equation}{0}

\section{Simulations}

We now demonstrate the performance of the proposed estimation
methods by simulations. We will compare them with some existing
methods including SIR (Li, 1991), SAVE (Cook and Weisberg, 1991),
PHD (Li, 1992) and rMAVE (Xia, Tong, Li and Zhu, 2002). The
computer codes used here can be obtained from
\textit{www.jstatsoft.org/} \textit{v07/i01/} for SIR, SAVE and
PhD methods (Courtesy of Professor S. Weisberg) and
\textit{www.stat.nus. edu/\~{}ycxia/} for rMAVE, dOPG and dMAVE.
In the following calculations, we use the quadratic kernel $ H(v)
= K_0(v^2) = (15/16) (1-v^2)^2 I(v^2 < 1)$ and $ \omega_0 = 0.01$.
The bandwidths in (\ref{band0}) and (\ref{band1}) are used. For
the inverse regression methods, the number of slices is chosen
between 5 to 30 that is most close to  $ n/(2p) $. We define an
overall estimation error of estimator $ \hat B: \hat B^\top \hat B
= I_q $ by the maximum singular value of $ B_0 B_0^\t - \hat B
\hat B^\t $; see Li et al (2004).

\begin{Example}\label{LIKC} \rm
Consider model
\beginy
  Y = \mbox{\rm sign}(2X^\t \beta_1 + \varepsilon_1)\log(|2X^\t \beta_2 + 4 +
  \varepsilon_2|),\label{model41}
\endy
where  $ \mbox{\rm sign}(\cdot) $ is the sign function.
Coordinates $ X \sim N(0, I_p)$, unobservable noises  $
\varepsilon_1 \sim N(0,1) $ and $\varepsilon_2 \sim N(0,1) $ are
independent. For $\beta_1 $, the first 4 elements  are all 0.5 and
the others are zero. For $\beta_2 $, the first 4 elements are
0.5,-0.5,0.5,-0.5 respectively and all the others are zero. A
similar model was investigated by Chen and Li (1998). In order to
show the effect on the estimation performances   of the number of
covariates, we vary  $ p $ in the simulation. With different
sample sizes, 200 replications are drawn from the model. The
calculation results are listed in Table 1. To get an intuition
about the quantity of estimation errors, Figure \ref{fli} shows a
typical sample of size $ n = 200$ and its estimate with estimation
error $ 0.21 $. The structure can be estimated quite well in the
sample.

\begin{figure}[ht]
\centerline{ \epsfxsize 5.5in \epsffile{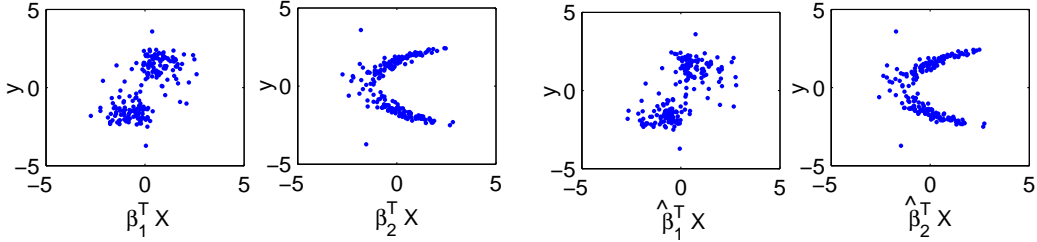}}
\vspace{-0.4cm} \caption{\small \it A typical data of size 200
from Example \ref{LIKC} with $ p = 10$ to show the quantity of
estimation error and its graphic performance. The left two panels
are plots of $ y $ against the true CS directions; the right two
panels  $ y $ against the estimated directions using dMAVE. The
estimated directions are respectively $ \hat\beta_1 = (0.42$,
0.64, 0.44, 0.45, -0.01, -0.07, 0.02, -0.00,   -0.08, 0.07$)^\t$
and $\hat\beta_2 = (-0.54,$    0.43,   -0.57,    0.43, 0.01,
-0.04, -0.01, 0.07, -0.05,
  0.07$)^\t$ with estimation error $ 0.21$.
 } \label{fli}
\end{figure}

\def\rt{\raisebox{-1.5ex}[0pt]}
{\footnotesize
\[
\begin{tabular}{|cc|cccccc|}
\multicolumn{8}{c}{\normalsize Table 1: Mean (and standard
deviation) of estimation errors for Example \ref{LIKC}}\\
\hline
$ n $ & $p$ & dOPG & dMAVE & rMAVE  & SIR & SAVE & PHD\\
\hline
     & 5 & \  0.25(0.09)   &  0.22(0.08)  &   0.43(0.19)  & 0.29(0.09)  &  0.87(0.19)   &   0.72(0.22)  \\
 100 & 10 & \ 0.55(0.19) &0.35(0.07) &   0.64(0.19)  & 0.46(0.10) & 0.94(0.06) & 0.90(0.13) \\
     & 20 & \ 0.81(0.13) &0.54(0.10) & 0.88(0.12)          & 0.64(0.11) & 0.96(0.06) & 0.93(0.07) \\ \hline
     & 5 &  \ 0.17(0.05) &0.14(0.04)  &  0.27(0.13)  & 0.19(0.05) & 0.55(0.26) & 0.47(0.15)\\
 200 & 10 & \ 0.32(0.09) &0.24(0.06) & 0.46(0.17)           &  0.30(0.06) & 0.96(0.08) & 0.73(0.16) \\
     & 20 & \  0.62(0.15) &0.36(0.06)  &  0.66(0.16)    &  0.43(0.06)& 0.93(0.04) &  0.94(0.08)\\ \hline
     & 5  & \ 0.13(0.04)     &0.13(0.04) &  0.19(0.07)   &  0.16(0.05) & 0.32(0.18) & 0.37(0.12)\\
 300 & 10 & \ 0.24(0.06)            &0.18(0.04) & 0.36(0.16)  & 0.24(0.05) & 0.85(0.17)  & 0.59(0.15) \\
     & 20 & \ 0.48(0.13) &0.28(0.05) & 0.55(0.16)  & 0.35(0.05) & 0.92(0.03)  & 0.84(0.12) \\ \hline
     & 5 &  \ 0.11(0.04)  &0.11(0.04) &  0.21(0.12)  & 0.14(0.04) & 0.22(0.11)  & 0.31(0.10)\\
 400 & 10 & \ 0.21(0.04) &0.16(0.04) & 0.31(0.11)  &0.21(0.05)  & 0.66(0.22) & 0.51(0.13) \\
     & 20 & \ 0.31(0.06) &0.25(0.04) &  0.49(0.15)  & 0.29(0.04) &  0.98(0.04) & 0.76(0.14) \\
\hline
\end{tabular}
\]
}

In model (\ref{model41}), the CS directions are hidden in a
complicated structure and are not easy to be detected directly by
the conditional regression mean function. When sample size is
large $(\ge 200)$ and $ p $ is not high ($=5$), all the methods
have accurate estimates. As $ p $ increases, rMAVE performs not so
well because the second direction is not captured in the
regression mean function; SAVE and PHD also fail to give accurate
estimates. SIR performs much better in all the situations than
SAVE and PHD. dOPG has about the same performance as SIR. dMAVE is
the best in all situations among all the methods.

\end{Example}

\begin{Example}\label{Emv} \rm Now, consider the CS in conditional mean
as well as the conditional variance as in the following model
\beginy
 && Y = 2(X^{\t}\beta_1)^d + 2\exp(X^\t\beta_2)\varepsilon,
  \label{mv}
\endy
where $ X = (x_1, \cdots, x_{10})^\t $ with $ x_1, \cdots, x_{10}
\sim Uniform(-\sqrt{3}, \sqrt{3}) $ and $ \varepsilon \sim N(0, 1)
$ are independent, $ \beta_1 = (1, 2, 0, 0, 0, 0, 0, 0, 0, 2)^\t/3
$ and $ \beta_2 = (0, 0, 3, 4, 0, 0, 0, 0, 0, 0)^\t$ $/5 $. For
model (\ref{mv}), one CS direction is contained in the regression
mean and the other in the conditional variance. One typical data
with size 200 is shown in Figure \ref{fmv}. Table 2 lists the
calculation results of 200 replications.

\begin{figure}[ht]
\centerline{ \epsfxsize 5.5in \epsffile{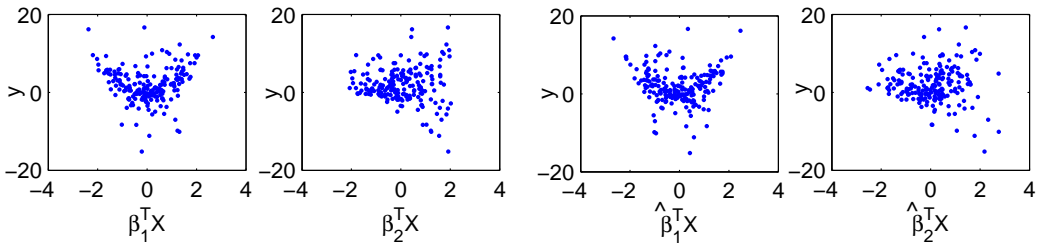} }
\vspace{-0.4cm} \caption{\small  \it A typical data with $ n = 200
$ from Example \ref{Emv} and its dMAVE estimation. The left two
panels are plots of $ y $ against the true CS directions
respectively; the right two panels $ y $ against the estimated
directions respectively with estimation error 0.31. } \label{fmv}
\end{figure}

Because rMAVE  cannot detect the CS directions  hidden in the
conditional variance directly, it has very poor overall estimation
performance as listed in Table 2. If $ d = 1$, i.e. the regression
mean function is monotonic, SIR works reasonably well; if $ d= 2
$, the regression mean function is symmetric and SIR fails to find
the direction hidden in the regression mean. As a consequence, its
performance is very poor. The performances of SAVE and PHD are
also far from satisfactory though they are applicable to the model
theoretically. The proposed dOPG and dMAVE perform very well and
are better than the existing methods listed in the table.

\def\rt{\raisebox{1.5ex}[0pt]}
\def\z{\hspace{-.3cm}}
{\footnotesize
\[
\begin{tabular}{|cc|cccccc|}
\multicolumn{8}{c}{\normalsize Table 2: Mean (and standard
deviation) of estimation errors for Example \ref{Emv} }\\
\hline
d & $ n $  & \ \ dOPG& dMAVE & rMAVE & SIR & SAVE & PHD\\
\hline
             & \ 100 &  \ \ 0.57(0.15)  &    0.44(0.12)  &   0.85(0.13) & 0.63(0.15) & 0.93(0.08) & 0.99(0.08)\\
 1         & \ 200 &  \ \ 0.36(0.08)  &   0.28(0.06)  &   0.76(0.16) & 0.42(0.09) & 0.91(0.12) & 0.98(0.07)\\
             & \ 400 &  \ \ 0.24(0.05) &    0.18(0.04)  &   0.68(0.15) & 0.29(0.06)& 0.64(0.16) & 0.97(0.07)\\
\hline
             &  \ 100 &  \ \ 0.63(0.19)  &  0.46(0.16)  &   0.85(0.16) & 0.96(0.09) & 0.90(0.06) & 0.91(0.11)\\
 2 &  \  200 &   \ \ 0.33(0.10) &   0.28(0.06)  &  0.70(0.18)& 0.95(0.07)  & 0.87(0.11) & 0.88(0.11) \\
            & \ 400 &  \ \  0.22(0.05) &
                0.19(0.04)  &   0.66(0.19) & 0.95(0.09) & 0.85(0.12) & 0.89(0.11)\\
\hline
\end{tabular}
\]
}

\end{Example}

\begin{Example}\label{rootn} \rm
In this example, we demonstrate the consistency rates of the
estimation methods by checking how the estimation errors change
with sample size $ n $. Consider model
\beginy
Y = \frac{x_1}{0.5 + (1.5+x_2)^2} + x_3(x_3 + x_4 + 1) + 0.1
\varepsilon, \label{model43}
\endy
where $ \varepsilon \sim N(0,1) $ and $ X \sim N(0, I_{10}) $ are
independent. Model (\ref{model43}) is a combination of the two
examples in Li (1991). For this model, all the theoretical
requirements for the methods are fulfilled. Therefore, it is fair
to use the model to check their consistency rates.

\begin{figure}[ht]
\centerline{ \epsfxsize 5.5in \epsffile{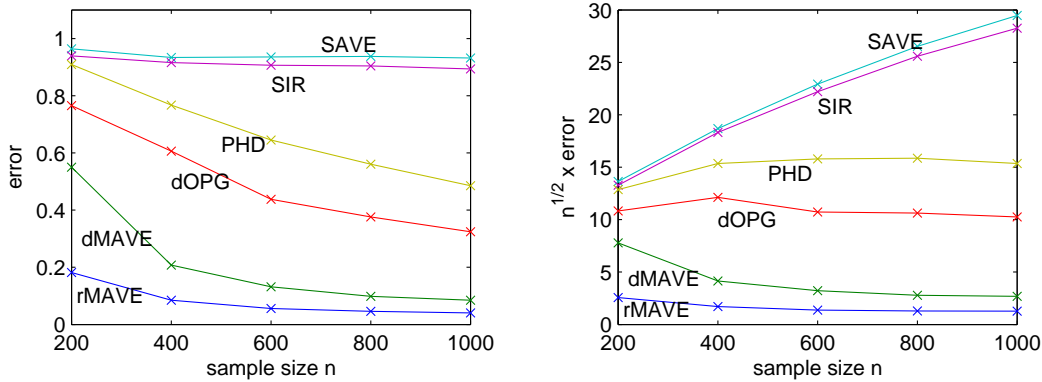} }
\caption{\small \it The calculation results for Example
\ref{rootn} using different estimation methods. The lines are the
mean of estimation errors with different sample size and 200
replications. The left panel is the plot of the errors against
sample size; the right panel is the errors multiplied by root-$n$
against sample size. } \label{figrootn}
\end{figure}

In the left panel in Figure \ref{figrootn}, the proposed methods
have much smaller estimation errors than the inverse regression
estimations. Because all the directions are hidden in the
regression mean function, it is not surprising that  rMAVE has the
best performance. Multiplied by root-$n$, the errors should keep
in a constant level if the theoretical root-$n$ consistency is
applicable to the range of sample size. The right panel suggests
that the estimation errors of SIR and SAVE do not start to show a
root-$n$ decreasing rate for the sample size up to 1000, while
PHD, rMAVE, dOPG and dMAVE demonstrate a clear root-$n$
consistency rate.

\end{Example}

\begin{Example}\label{circle} \rm In our last example, we consider a model with
a very complicated structure. Suppose $ (X_i, Y_i), i = 1, 2,
\cdots, n, $ are drawn independently from model $Y = \beta_1^\top
X/2 +\varepsilon(1-|\beta_1^\top X|^2)^{1/2} $, where $ (X,
\varepsilon)$ satisfies $ \{ X\sim N(0, I_{10}),  \varepsilon \sim
N(0, 1): |\beta_1^\top X|\le 1, |\beta_2^\top X|\le 1, 0.5 <
(\beta_1^\top X)^2 (1-\varepsilon^2) + \varepsilon^2 \le 1 \} $, $
\beta_1 $ and $ \beta_2 $ are defined in Example \ref{LIKC}. The
calculation results based on 200 replications are listed in Table
3. Because of the complicated structure as shown in Figure
\ref{complex}, the CS directions are not easy to be estimated and
observed directly. However, with moderate sample size, the
proposed methods can still estimate the directions accurately. It
is interesting to see that SAVE also works in  this example.

\

{\footnotesize
\[
\begin{tabular}{|c|cccccc|}
\multicolumn{7}{c}{\normalsize Table 3: Mean (and standard
deviation) of estimation errors for Example \ref{circle}}\\
\hline
 $ n $  & dOPG& dMAVE & rMAVE & SIR & SAVE & PHD\\
\hline
 200 & 0.5909(0.29)  &  0.5089(0.30)  &  0.9411(0.07)  &  0.8770(0.12) &   0.9242(0.19)  & 0.9833(0.05)  \\
 400 &  0.2117(0.19) &   0.1498(0.10) &   0.9573(0.05) &   0.8783(0.13) & 0.7677(0.18) &   0.9789(0.03) \\
 600 &  0.1148(0.04) &   0.1059(0.03)  &  0.9725(0.03)  &  0.8758(0.13) &
0.5357(0.21) &   0.9799(0.03)  \\
 800 &   0.0876(0.03) &   0.0862(0.02) &   0.9744(0.03)  &  0.8737(0.14) &
0.3657(0.13) &   0.9757(0.04) \\
1000&   0.0782(0.02)  &  0.0779(0.02)  &  0.9671(0.04)  &
0.8819(0.13) &  0.2604(0.06) & 0.9789(0.04)
\\
\hline
\end{tabular}
\]
}

Based on the simulations, we have the following observations. (1)
The existing methods (rMAVE, PHD, SIR and SAVE) fail in one way or
another to estimate the CS directions efficiently, while dOPG and
dMAVE are efficient for all the examples. (2)  dOPG and dMAVE
demonstrate very good finite sample performance, even a root-$n$
rate of estimation efficiency, while some of the existing methods
do not show a clear root-$n$ rate in the range of sample sizes
investigated. (3) dOPG and dMAVE are less sensitive to the number
of covariates than PHD, SAVE and SIR. Simulations not reported
here also suggest that the asymmetric design of $ X$ has less
effect on dOPG and dMAVE than that on the inverse regression
estimations. (4) If the CS  directions are all hidden in the
regression mean function, rMAVE is the best and should be used.
Otherwise, dOPG and dMAVE are recommended.

\begin{figure}[ht]
\centerline{ \epsfxsize 5in \epsffile{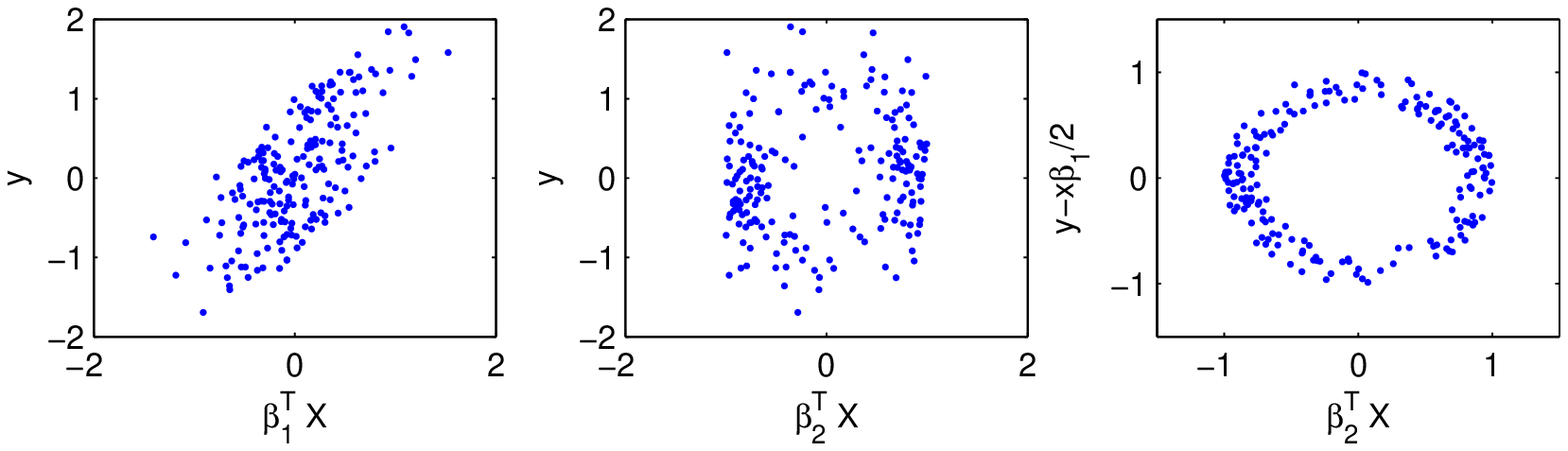} }
 \centerline{ \epsfxsize 5in \epsffile{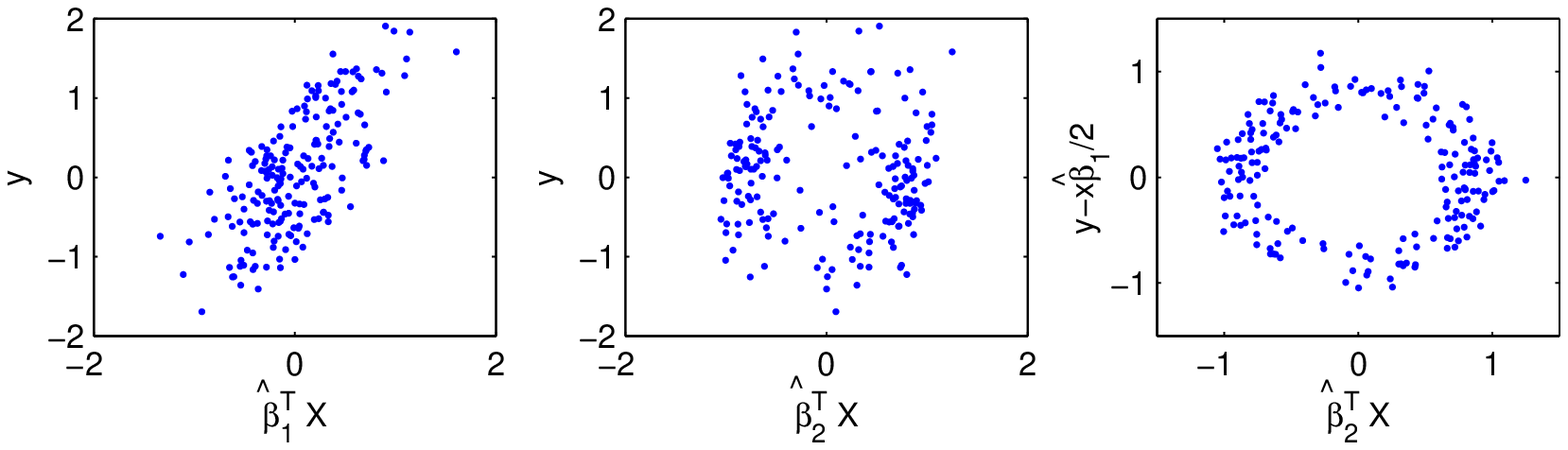} }
\vspace{-0.3cm} \caption{\small \it A typical data from Example
\ref{circle} with $ n = 200 $ and its dMAVE estimation. The upper
three panels are plots of $ y $ against the true CS directions and
$ y - x^\top \beta_1/2$ against the second direction respectively;
the lower three panels are plots of $ y $ against the estimated CS
directions (with estimation error 0.32) and $ y - x^\top \hat
\beta_1/2$ against the  second estimated direction respectively. }
\label{complex}
\end{figure}

\end{Example}

\section{Real data analysis}

\begin{Example}[Cars data]\label{cars} \rm
This data was used by the American Statistical Association in its
second (1983) exposition of statistical graphics technology. The
data set is available at {\footnotesize
http://lib.stat.cmu.edu/datasets/cars.data}. There are 406
observations on  8 variables: miles per gallon ($Y$), number of
cylinders ($X_1$), engine displacement ($X_2$), horsepower
($X_3$), vehicle weight ($X_4$), time to accelerate from 0 to 60
mph ($X_5$), model year ($X_6$), and origin of a  car (1.
American, 2. European, 3. Japanese).

Now we investigate the relation between response variable $Y$ and
covariates $ X = (X_1, \cdots, X_8)^\top $, where $ X_1, \cdots,
X_6 $ are defined above,  $ X_7=1 $ if a car is from America and 0
otherwise; $ X_8=1$ if it is from Europe and 0 otherwise. Thus,
$(X_7, X_8)=(1,0), (0,1)$ and (0,0) correspond to American cars,
European cars and Japanese cars respectively. For ease of
explanation, all covariates are standardized separately.  When
applying dOPG to the data, the first 4 largest eigenvalues are
21.1573, 1.6077, 0.2791 and 0.2447 respectively. Thus, we consider
CS with dimension 2. Based on dMAVE, the two directions
(coefficients of $ X$) are estimated as $ \hat\beta_1 = ($-0.33,
-0.45, -0.45, -0.53, 0.14, 0.42, 0.00, -0.02$)^\top $ and $
\hat\beta_2 = ($0.00, 0.15, -0.10, -0.23, -0.12, -0.17, -0.88,
0.29$)^\top $ respectively. The plots of $ Y $ against $ \hat
\beta_1^\top X $ and $ \hat \beta_2^\top X $ are shown in Figure
\ref{figcars}.

\begin{figure}[ht]
\centerline{ \epsfxsize 6.5in \epsffile{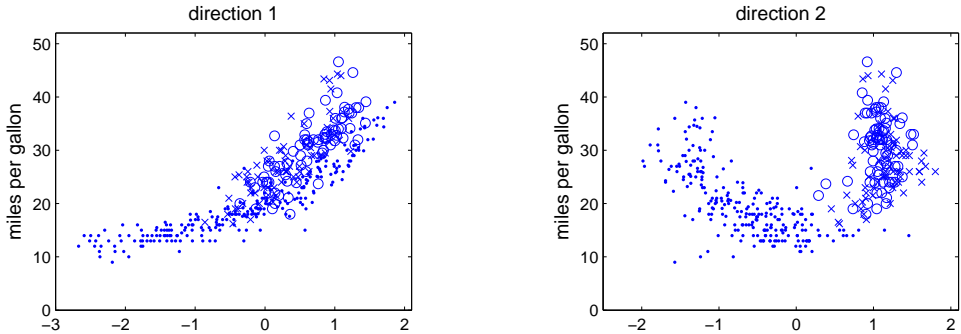} }
\vspace{-5.9cm} \caption{\small  \it The estimation results for
Example \ref{cars} using dMAVE. The two panels are plots of $ Y $
against the two estimated CS directions respectively.  The origins
of cars are denoted by ``$\cdot$'' for American cars, ``$\times$''
for European cars, and ``$\circ$'' for Japanese cars. }
 \label{figcars}
\end{figure}

Based on the estimated CS directions and Figure \ref{figcars}, we
have the following insights to the data. The first direction
highlights the common structure for cars of all origins: miles per
gallon ($Y$) decreases with number of cylinders ($X_1$), engine
displacement ($X_2$), horsepower ($X_3$) and vehicle weight
($X_4$),  and increases with the time to accelerate ($X_5$) and
model year ($X_6$). The second direction indicates the difference
between American cars and European or Japanese cars.

\end{Example}

\begin{Example}[Ground level Ozone]\label{ozone} \rm
Air pollution has serious impact on the health of plants and
animals (including humans); see the report of the World Health
Organization (WHO) (2003). Substances not naturally found in the
air or at greater concentrations than usual are referred to as
``pollutants''. The main pollutants include nitrogen dioxide
(NO$_2$), carbon dioxide (CO), sulphur dioxide (SO$_2$),
respirable particulates, ground-level ozone (O$_3$) and others.
Pollutants can be classified as either primary pollutants or
secondary pollutants. Primary pollutants are substances directly
produced by a process, such as ash from a volcanic eruption or the
carbon monoxide gas from a motor vehicle exhaust. Secondary
pollutants are products of reactions among primary pollutants and
other gases. They are not directly emitted and thus cannot be
controlled directly. The main secondary pollutant is ozone.

Next, we investigate the statistical relation between the level of
ground-level ozone with the levels of primary pollutants and
weather conditions by applying our method to  the pollution data
observed in Hong Kong (1994-1997, {\footnotesize http://www.hku.hk
/statistics/paper/}) and Chicago (1995-2000, {\footnotesize
http://www.ihapss.jhsph.edu/data/data.htm}). This investigation is
of interest in understanding how the secondary pollutant ozone is
generated from the primary pollutants and weather conditions. Let
$ Y $, $ N, S, P, T $ and $ H $ be the weekly average levels of
ozone, nitrogen dioxide (NO$_2$), sulphur dioxide (SO$_2$),
respirable particulates, temperature and humidity respectively. To
include the interaction between  primary pollutants and weather
conditions into  the model directly, we further consider their
cross-products resulting in 15 covariates all together, denoted by
$ X$. For ease of explanation, all covariates are standardized
separately. For all possible working dimensions, only the first
two dimensions show clear relations with $ Y$. We further
calculate the eigenvalues in dOPG. The largest four eigenvalues
are $ 10.78, 2.93, 2.11, 1.70$ respectively for Chicago, and $
6.89, 1.24,    0.69, 0.52 $ for Hong Kong. Now we consider the
dimension reduction with efficient dimension 2 although the
estimation of the number of dimension needs further investigation.
The estimates for the first two directions are given in Table 4.

\def\rt{\raisebox{0.5ex}[0pt]}
{\footnotesize
\[
\begin{tabular}{|c|c|cccccccc|}
\multicolumn{10}{c}{\normalsize Table 4: The estimated CS directions in Example \ref{ozone}}\\
\hline
 City & Direction & $N$ & $S$ &  $P$ &  $T$
   &  $H$ &  $N*S$&  $N*P$ &  $N*T$ \\ \hline
  & $\beta_1$ &   0.10&
   -0.13 &
   -0.06 &
   -0.00 &
   -0.00 &
    0.06 &
    0.29 &
    0.19
\\
\rt{Chicago}  & $\beta_2$ & -0.10 &
   -0.11 &
    0.39 &
   -0.25 &
   -0.07 &
    0.12 &
   -0.15 &
    0.09
   \\
 & $\beta_1$ &     0.32 &
   -0.15 &
    0.23 &
    0.10 &
   -0.41 &
   -0.07 &
    0.20 &
    0.42
   \\
\rt{Hong Kong}   & $\beta_2$ &     -0.04 &
   -0.08 &
   -0.12 &
    0.18 &
    0.19 &
   -0.21 &
    0.35 &
    0.17
\\
\hline\hline city  & Direction &  $N*H$ &  $S*P$ & $S*T$    &
$S*H$ & $P*T$ &  $P*H$
   &  $T*H$ & \\ \hline
  & $\beta_1$ &        0.04 &
   -0.18 &
    0.27 &
   -0.01 &
   -0.06 &
    0.36 &
    0.77 &
 \\
\rt{Chicago}  & $\beta_2$ &    -0.51 &
    0.46 &
   -0.20 &
   -0.21 &
   -0.15 &
   -0.16 &
    0.32 &
\\
  & $\beta_1$  &      0.10 &
    0.01 &
   -0.05 &
   -0.31 &
    0.53 &
    0.12 &
   -0.14 &
  \\
\rt{Hong Kong}  & $\beta_2$  &
   -0.52 &
   -0.26 &
   -0.18 &
    0.42 &
    0.22 &
   -0.29 &
   -0.19 &
 \\ \hline
\end{tabular}
\]
}

The  plots of $ Y $ against the two estimated directions are shown
in Figure \ref{figozone}. The plots show strong similar patterns
in the two separated cities. If we check the estimated
coefficients (directions), NO$_2$ and particulates (or their
interaction) are the most important pollutants that affect the
level of ozone. Temperature and humidity and their interaction are
the other important factors. The interactions of weather
conditions with NO$_2$ and particulates also contribute to the
variation of ozone levels. These statistical conclusions give
support to  the chemical claim that ozone is formed by chemical
reactions between reactive organic gases and oxides of nitrogen in
the presence of sunlight; see the report of WHO (2003).

\begin{figure}[ht]
 \centerline{ \raisebox{15ex}[0pt]{Hong Kong} \ \epsfxsize 4in \epsffile{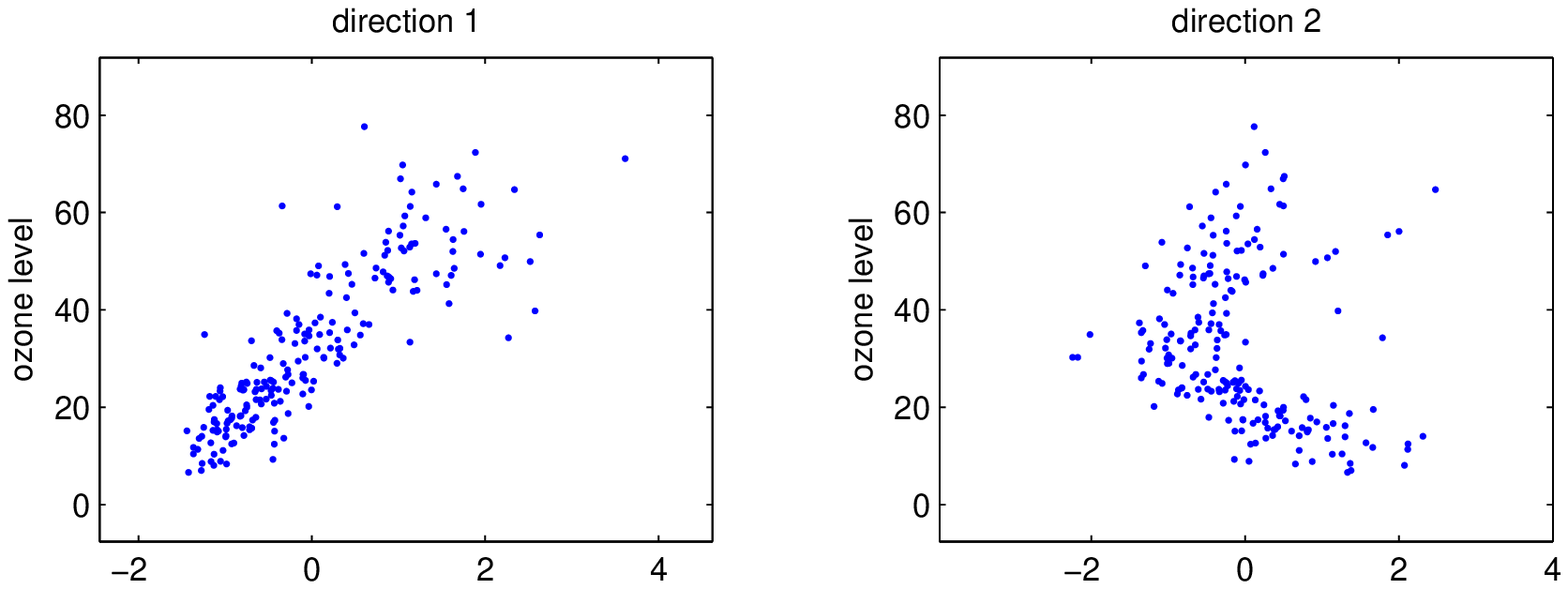} }
 \centerline{ \raisebox{15ex}[0pt]{Chicago \ \ \ } \ \epsfxsize 4in \epsffile{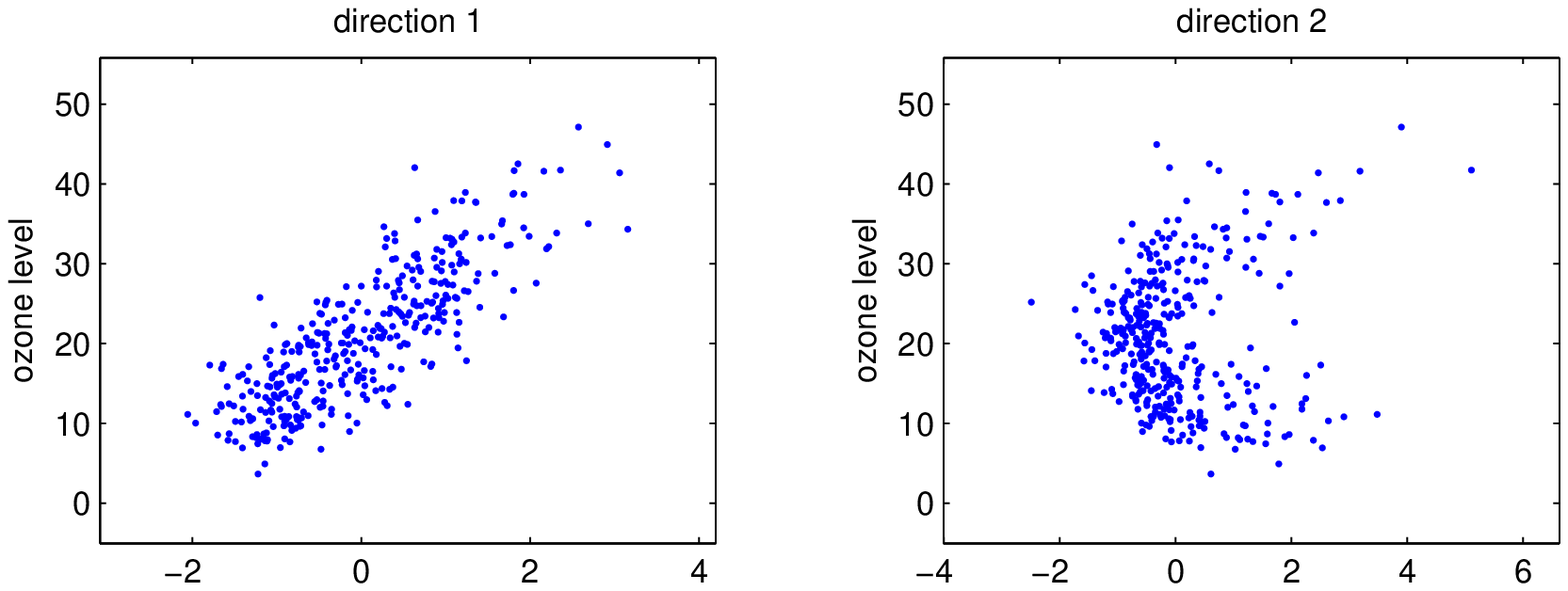} }
\caption{\small \it The estimation results for Example \ref{ozone}
using dMAVE. The upper two panels are the levels of ozone against
the first two estimated CS directions in Hong Kong, the lower two
panels are those in Chicago. } \label{figozone}
\end{figure}

\end{Example}




\vspace{-0cm}

\setcounter{equation}{0}
\section{Proofs}
\subsection{Basic ideas of the proofs}

The basic idea to prove the theorems is based on the convergence
of the algorithms and that the true dimension reduction space is
the attractor of the algorithms. We here give a more detailed
outline for the proof of Theorem 3.2. Suppose the estimate of $
B_0 $ in an iteration of the dMAVE algorithm is $ B_{(t)} $. It
follows from Step 2 that
\beginy
\n \b^{(t+1)} &=& \ell(B_0) + \Big\{\sum_{k,j,i=1}^n
\rho^{(t)}_{jk} K_{h_{t}}( B_{(t)}^{\t} X_{ij})
        X^{(t)}_{ijk} (X^{(t)}_{ijk})^{\t}
       \Big\}^{-1}\\
       && \hspace{1cm}\times \sum_{k,j,i=1}^n  \rho^{(t)}_{jk} K_{h_{t}}(
B_{(t)}^{\t} X_{ij})
        X^{(t)}_{ijk}\{H_{b,i}(Y_k) - a^{(t)}_{jk} - \ell( B_0)^\t
        X^{(t)}_{ijk}\},  \qquad\label{relation0}
\endy
where $ X_{ijk}^{(t)} $ is defined in the algorithm. By the
decomposition in Step 3, we obtain estimate $ B_{(t+1)} $ in the
next iteration. If the initial value $ B_{(1)}$ is a consistent
estimator of $ B_0 $, by Lemmas \ref{krq}, \ref{denominator} and
\ref{numerator} below, we will obtain a recurring relation for the
iterations as
\beginn
\ell(B_{(t+1)}) - \ell( B_0) = \Theta_t \{ \ell(B_{(t)}) - \ell(
B_0)\} + \Gamma_{n,t}, \label{relation}
\endn
with $  | \Theta_t | < 1 $ and $  |\Gamma_{n,t}| = o(1) $ almost
surely when $ t \ge 1 $. Therefore, the dimension reduction space
is an attractor in the algorithm. This recurring relation is then
used to prove the convergence of the algorithm and  the
consistency of the final estimator.   To ensure the convergence of
the algorithm, we need to consider the consistency with
probability 1.

The details of the proofs are organized as follows. In section
6.2, we first list a series of lemmas, Lemmas 6.1-6.5. Based on
these Lemmas the theorems are then proved. The proofs of Lemmas
6.1-6.5 are algebraic albeit complex calculations of Lemmas 6.6
and 6.7. They can be found in Xia (2006b) and are available upon
request. Lemmas 6.6 and 6.7 are two basic results used in the
proof dealing with the uniform consistency. Their proofs are given
in section 6.3.

\subsection{Proofs of the theorems}

We first introduce a set of notations. Let $ \varepsilon_{b, i}(y)
= H_b(Y_i-y) - E(H_b(Y_i-y)| X_i) $,  $ \D_{_Y} \subset \R $ be a
compact interior support of $ Y $, i.e. for any $ v \in \D_{_Y} $,
there exists $ \delta> 0 $ such that $ \inf_{y: |y-v|<\delta}
f_{_Y}(y) > 0 $. Similarly, we can define a compact interior
support  $ \D_{_X} $ for $ X $. For $ {{\cal B}} \subset  \{ B:
B^\top B = I_q\} $,  define $ \delta_B = \max\{ |B - B_0|: B \in
{{\cal B}}\}$. For any index set $ {{\cal Z}} $ and random matrix
$ A_n(z) $, we say $A_n(z)= \O(a_n|z \in {{\cal Z}}) $, or $
A_n(z) = \O(a_n) $ for simplicity, if $ \sup_{z \in {{\cal Z}}}
|A_n(z)| / a_n = O(1)$ almost surely. As usual, $ A_n= O_P(a_n) $
indicates that every term in $ A_n $ is $ O(a_n) $ in probability
as $ n \to \infty $. Recall that $ B_0 = (\beta_{01}, \beta_{02},
\cdots, \beta_{0q})$ and $ B = (\beta_{1}, \beta_{2}, \cdots,
\beta_{q})$. Let $ H_{b,i}^{1,B}(x) \ = \ g_b(B_0^\t x, y) + \btdt
g_b( B_0^\t x, y) B^\t X_{ix}$, \ $ H_{b, i}^{2,B}(x) \ = \
\sum_{\iota,\kappa=1}^q \btd^2_{\iota,\kappa}
   g_b( B_0^\t x, y)(\beta_\iota^\t X_{ix}) (\beta_\kappa^\t
   X_{ix})/2$ and $ H_{b, i}^{3,B}(x) =  \sum_{\iota,\kappa, \tau=1}^q
\{\btd^3_{\iota,\kappa,\tau}
   g_b(B_0^\t x, y)$ $(\beta_\iota^\t X_{ix}) (\beta_\kappa^\t
X_{ix}) (\beta_\tau^\t X_{ix})\}/6$, where $ X_{ix} = X_i - x $, $
\btd g_b(v_1, $ $\cdots, v_q, y) $ is defined in Section 2 and
$$
\btd^2_{\iota,\kappa}  g_b(v_1, \cdots, v_q, y) = \frac{\partial^2
}{\partial v_\iota \partial v_\kappa} g_b(v_1, \cdots, v_q, y)
\dfor \iota,  \kappa =1,2,\cdots, q,
$$
and  $ \btd_{\kappa, \tau, \iota}^3 g_b $ is defined naturally. By
Taylor expansion of $ g_b( B_0^\t X_i, y) $ at $ B_0^\t x $, it
follows from  model (\ref{model}) that
\beginy
  H_{b, i}(y)  = H_{b,i}^{1,B_0}(x) +
  H_{b, i}^{2,B_0}(x)+ H_{ b, i}^{3,B_0}(x)+
  \varepsilon_{b, i}(y)   + O( |B_0^\t X_{ix}|^4) \qquad  \label{Hx}
\endy
almost surely. Let $ \delta_{mh} = (nh^m/\log n)^{-1/2} $, $
\delta_{mhb} = (nh^m b/\log n)^{-1/2} $ for any integer $ m $, $
\delta_b = (nb/\log n)^{-1/2}$, $ \delta_n = (\log n/n)^{1/2}$ and
$ r_{mhb} = h^2 + b^4 + \delta_b + \delta_{mh} $. Let $ f_B, f $
and $ f_{_Y} $ be the density functions of $ B^\t X $, $ X$ and $
Y$ respectively. Again, for simplicity, we write $ f_B(x), \mu_B(
x), w_B(x) $ for $ f_B(B^\top x),  \mu_B(B^\top x) $ and $
w_B(B^\top x) $ respectively; see also the definitions in Section
3. Let $ c, c_0, c_1, \cdots, $ be a sequences of positive
constants, while $ c$ may have different values at different
places.

\begin{Lemma}\label{krp} \rm [Kernel smoother in the first iteration] Let
\beginn
{a_{xy}\choose b_{xy}h} = \Big\{\sum\limits_{i=1}^n
   K_{h}(X_{ix})
   {1 \choose  X_{ix}/h } {1 \choose  X_{ix}/h }^{\t}\Big\}^{-1}
   \sum\limits_{i=1}^n K_{h}(X_{ix})
   {1 \choose  X_{ix}/h } H_{b,i}(y). \label{exp1}
\endn
Under assumptions (C1), (C2) and (C4), if $ h \to 0, b \to 0 $ and
$ n h^{p+2} b /\log n \to \infty $, then we have
\beginn
 a_{xy} \eqc g_b(B_0^\top x, y) + \frac12 \sum\limits_{\kappa = 1}^q \btd_{\kappa, \kappa}^2 g_b(B_0^\top x, y)
    h^2 + \O(h^3 + \delta_{phb}|x\in \D_{_X}, y\in  \D_{_Y}),  \\
 b_{xy} \eqc B_0\btd g_b(B_0^\t x,y) + \{\mu_{2p} nh^2 f(x)\}^{-1}
\sum\limits_{i=1}^n
  K_{h}(X_{ix}) X_{ix} \varepsilon_{b,i}(y) \\
  &&\hspace{6cm}+
  \O(h^2+\delta_{phb}|x\in \D_{_X}, y\in  \D_{_Y}).
\endn
\end{Lemma}

\begin{Lemma} \label{krqp} \rm [Kernel smoother in dOPG] Define ${{\cal D}}_q = \{D =  B
\diag(\lambda_{1}, \cdots,$ $ \lambda_{q}) B^\t + \tilde B
\diag(\lambda_{q+1}, $ $\cdots, \lambda_{p})\tilde B^\top$:  $ (B,
\tilde B)^\t (B, \tilde B) = I_p $, $c_1> \min(\lambda_1, \cdots,
\lambda_q ) \ge c_0 >0 $, $ B \in {{\cal B}} $ and $ \max(
\lambda_{q+1}, $ $ \cdots, \lambda_p )/h^2 \le e_n \} $. Let
\beginn
 S_n^{D}(x) = n^{-1}\sum_{i=1}^n
   K_{h}( D^{1/2} X_{ix})
   {1 \choose   X_{ix} } {1 \choose   X_{ix} }^{\t}
\endn
and
\beginn
&& \hspace{-0.8cm} {a_{xy}^{D} \choose b^{D}_{xy}} = \{n S_n^ {D}
(x)\}^{-1}
   \sum_{i=1}^n  K_{h}( D^{1/2} X_{ix})
   {1 \choose  X_{ix}} H_{b,i}(y).
\endn
Under assumptions (C1), (C2) and (C4), if $ nh^{q+2} b/\log n \to
\infty $,  $b \to 0, h\to 0$, $ \delta_B/h \to 0 $ and $ e_n \to
0, $ then we have
\beginn
   a^{D}_{xy} \eqc  g_b(B_0^\top x, y) + \frac12 \sum\limits_{\kappa = 1}^q \btd_{\kappa, \kappa}^2 g_b(B_0^\top x, y)
    h^2 + \O(h^3 + \delta_{qhb}|x\in \D_{_X}, y\in  \D_{_Y}, D\in {{\cal D}}_q), \ \ \ \\
   b^{D}_{xy} \eqc  B_0 \{  \btd  g_b( B_0^\t x, y)+ \O(h^2 +
\delta_{qh} + e_n)\}
      + \E^D_{n,0}(x,y)  \\
      && \hspace{6.8cm}+  \O( \epsilon_{qhb}|x\in \D_{_X}, y\in  \D_{_Y}, D\in {{\cal D}}_q),
\endn
where $ \epsilon_{qhb}  = h^4 + (h^2  + \delta_{qh})\delta_{qhb} +
(h^2+\delta_{qhb})e_n  + (h + \delta_{qhb}/h)\delta_B $ and
\beginn
\E^D_{n,0}(x,y) = h^{p-q} \{ n f_B( x) \}^{-1} \prod_{\tau=1}^q
   \lambda_\tau^{1/2}
    \bar w_B^{+}(x)
    \sum_{i=1}^n K_{h}( D^{1/2}    X_{ix})
   \{\mu_B(x)-X_i\} \varepsilon_{b,i}(y).
\endn
\end{Lemma}


\begin{Lemma}\label{krq} \rm [Kernel smoother in dMAVE] Let
\beginn
 \Sigma_n^B(x) = n^{-1}\sum_{i=1}^n
   K_{h}(B^\t X_{ix})
   {1 \choose  B^\t X_{ix}/h } {1 \choose  B^\t X_{ix}/h }^{\t}
\endn
and
\beginn
&& \hspace{-0.8cm} {a_{xy}^B\choose d^B_{xy}h} = \{n
\Sigma_n^B(x)\}^{-1}
   \sum_{i=1}^n  K_{h}( B^\t X_{ix})
   {1 \choose B^\t X_{ix}/h } H_{b,i}(y).
\endn
Under assumptions (C1), (C2) and (C4), if $ nh^q b/\log n \to
\infty $,  $b \to 0, h\to 0$ and $ \delta_B/h \to 0 $, then
\beginn
a^B_{xy} \eqc g_b( B_0^\t x, y) + \btdt  g_b( B_0^\t x,
y)(B_0-B)^\t \nu_{_B}( x)  +
   \frac12 \sum_{\kappa=1}^q
   \btd^2_{\kappa,\kappa}g_b( B_0^\t x, y) h^2\\
   && \hspace{1.5cm} +\V_{1n}^B(x,y)
     + \O(h^4 + \delta_{qh} \delta_{qhb} + h\delta_B + \delta_B^2|x\in \D_{_X}, y\in  \D_{_Y}, B\in {{\cal B}}),\\
   d^B_{xy}h \eqc \btd  g_b( B_0^\t x, y) h + M_{1n}^B(x,y) h^3\hspace{-0.1cm} +
   \V_{2n}^B(x,y) \\
 && \hspace{3.5cm}+ \O(h^4\hspace{-0.15cm} + \delta_{qh} \delta_{qhb}
\hspace{-0.05cm}+ h\delta_B\hspace{-0.05cm} + \delta_B^2|x\in
\D_{_X}, y\in  \D_{_Y}, B\in {{\cal B}}),
\endn
where
$$ \V_{1n}^B(x,y) =   \{1+M_{2n}^B(x, h) h \}\E^B_{n,1}(x,y)
+ M_{3n}^B(x, h) h \E^B_{n,2}(x,y),
$$
$$
\V_{2n}^B(x,y) =  M_{4n}^B(x) h \E^B_{n,1}(x,y)  + \{1+
M_{5n}^B(x, h) h\} \E^B_{n,2}(x,y),
$$
$ M_{kn}^B(x), k = 1,2,
\cdots, 5, $ are bounded continuous functions (details can be
found in the proofs) and
\beginn
\E^B_{n,1}(x,y) =\{ n f_B( x)\} ^{-1} \sum_{i=1}^n K_h(B^\t
X_{ix}) \varepsilon_{b, i}(y),\hspace{0.5cm} \\
 \E^B_{n,2}(x,y) = \{ n h
f_B( x)\} ^{-1} \sum_{i=1}^n K_h(B^\t X_{ix}) B^\t X_{ix}
\varepsilon_{b, i}(y).
\endn
\end{Lemma}

\begin{Lemma}\label{denominator} \rm [Denominator of dMAVE] Let  $ \hat \rho^B_{jk} = \rho(\hat
f_B (X_j))\rho( \hat f_{_Y}(Y_k) ) $, where
$$
\hat f_B( x) = n^{-1} \sum_{i=1}^n K_h(B^\t X_{ix}), \qquad \hat
f_{_Y}(y) = n^{-1} \sum_{i=1}^n H_b(Y_i-y).
$$
Let $ X^{B}_{ijk} = d_{jk}^B \otimes X_{ij} $ where $ d_{jk}^B =
d_{X_j Y_k}^B $. Suppose (C1)--(C4) hold and $ nh^{q+2}b /\log n $
$\to \infty $,  $ nb^2/\log n \to \infty$, $b \to 0, h\to 0$ and $
\delta_B/h \to 0 $. We have
\beginn
\Big\{ n^{-3} \sum_{k,j,i=1}^n \hat \rho^B_{jk}
            K_h(B^{\t} X_{ij})
        X^{B}_{ijk} (X^{B}_{ijk})^{\t}\Big\}^{-1} \eqc \ (I_q\otimes B)
        L_1^B (I_q\otimes B^\top) h^{-2} + (I_q\otimes B)
        L_2\\
        && \hspace{-2.3cm} +
        L_3 (I_q\otimes B^\top) +  \frac12 D_B^+
         + \O\{(r_{qhb}+\delta_{qhb})/h| B\in {{\cal B}}),
\endn
where $ L_1, L_2 $ and $ L_3 $ are constant matrices (details can
be found in the proof) and  $ D_B= \int \rho(f_B( x))
 \rho(f_{_Y}(y)) \btd g_b( B_0^\t x, y) \btdt g_b( B_0^\t x, y)
\otimes \{\nu_{_B}(x) \nu_{_B}^\t (x)\}f(x) f(y) dx dy $.
\end{Lemma}

\begin{Lemma}\label{numerator} \rm [Numerator of dMAVE]
Suppose conditions  (C1)--(C4) hold. If $b \to 0, h\to 0$, $ nh^q
b/\log n \to \infty $,   $ nb^2/\log n \to \infty$ and $
\delta_B/h \to 0 $, then
\beginn
 n^{-3}\hspace{-0.2cm}\sum_{k,j,i=1}^n \hat \rho^B_{jk} K_h( B^{\t} X_{ij})
        X^{B}_{ijk}\{H_{b,i}(Y_k) - a^{B}_{jk} - \ell( B_0)^\t
        X^{B}_{ijk}\}= D_B ( \ell(B) - \ell(B_0)) \qquad \\
    \hspace{-5cm}+ \Phi_n(B_0) + \O\{h^4 + r_{qhb}
   \delta_{qhb}+ \delta_{qhb}^2 + \delta_n^2/b^2 +
   (\delta_{qhb}/h+h)\delta_B|B \in {{\cal B}}\},
\endn
where $ a^B_{jk} = a^B_{_{X_j Y_k}},  \Phi_n(B_0) =  O(
\delta_n+\delta^2_{qhb}/h) $ almost surely and $ \Phi_n(B_0) =
O_P(n^{-1/2}) $  with $ (I_q\otimes B_0^\top) \Phi_n(B_0) = 0 $
and $ \sqrt{n} \Phi_n(B_0) \stackrel{D}\to N(0, \Sigma_0), $ where
$ \Sigma_0  $ is given in Theorem \ref{mainMAVE}.
\end{Lemma}


{\bf Proof of Theorem \ref{mainOPG}} By Lemma \ref{krp}, write
 \beginn
b_{xy} = B_0 c_n(x, y)   + \{\mu_{2p} nh_0^2 f(x) \} ^{-1}
\sum_{i=1}^n K_{h_{_0}}(X_{ix}) X_{ix} \varepsilon_{b_{_0},i}(y) +
\tilde B_0 \O(h_0^2+\delta_{ph_0 b_0}),
 \endn
where $ (B_0, \tilde B_0) $ is a $p\times p $ orthogonal matrix
and $ c_n(x,y) = \btd g_b(B_0^\top x, y) + \O(h_0^2 +
\delta_{ph_{_0}b_{_0}}) $. By Lemma \ref{basic1}, the second term
on the right hand side above is $ \O(\delta_{ph_{_0}b_{_0}}/h_0)
$. It follows from step 2 in the dOPG algorithm that
\beginy
\n \hat \Sigma_{(1)} &=& (B_0, \tilde B_0) C_n (B_0, \tilde
B_0)^\top\hspace{-0.1cm} + n^{-3} \hspace{-0.2cm}\sum_{i,j,k=1}^n
( S_{ijk} + S_{ijk}^\top )\\
&& \hspace{5.5cm} + \O\{(h_0^2+\delta_{ph_0
b_0})\delta_{ph_{_0}b_{_0}}/h_0\}, \label{opga}
\endy
where $ \hat \Sigma_{(1)} \ $ and $ \ \rho_{jk}^{(0)}$ are defined
in the algorithm, $ S_{ijk} = \rho_{jk}^{(0)} \{\mu_{2p} h_0^2 f(
X_j)\}^{-1} $ $B_0 \btd g_{b_{_0}}( B_0^\t X_j, Y_k)
K_{h_{_0}}(X_{ij}) X_{ij}^\t \varepsilon_{b_{_0},i}(Y_k)$ and
\beginn
 C_n &=& n^{-2}\hspace{-0.1cm}\sum_{j,k=1}^n \rho_{jk}^{(0)}
{c_n (X_j, Y_k) \choose \O(h_0^2+\delta_{ph_0 b_0})}{ c_n(X_j,
Y_k) \choose \O(h_0^2+\delta_{ph_0 b_0})}^\top\hspace{-0.1cm} \\
 &=&
\left(\hspace{-0.2cm}
\begin{array}{cc} \Lambda_n^{(1)} & \O(h_0^2 + \delta_{ph_0b_0}) \\
\O(h_0^2 + \delta_{ph_0b_0}) & \O(h_0^4 + \delta_{ph_0b_0}^2)
\end{array}\hspace{-0.2cm} \right),
\endn
where $ \Lambda_n^{(1)} = n^{-2}\sum_{j,k=1}^n \rho_{jk}^{(0)} c_n
(X_j, Y_k)  c^\top_n(X_j, Y_k) $. By Lemma \ref{basic1}, we have $
 \tilde f_{_Y}^{(0)}(y) = f_{_Y}(y) + f_{_Y}''(y)b_0^2/2 +\O(b_0^4+ \delta_{b_{_0}}|y \in \D_{_Y}),
\tilde f^{(0)}(x) = f(x) +   \O(h_0^2+ \delta_{ph_{_0}}|x\in
\D_{_X})$. By the definition of $ \rho(.)$,  we have $
\rho_{xy}^{(0)} = \rho(f( x))\tilde \rho_{b_0}( f_{_Y}(y)) +
\O(r_{p h_0 b_0}|x \in {{\Bbb R}}^p, y \in {\Bbb R}) $, where $
\tilde \rho_{b_0}(f_{_Y}(y)) = \rho(f_{_Y}(y)) +\rho'(f_{_Y}(y))
f_{_Y}''(y)b_0^2/2 $. Let
\beginn
\tilde S_{ijk} &=& \rho(f(X_j))\tilde \rho_{b_0}( f_{_Y}(Y_k)) B_0
\btd g_{b_{_0}}( B_0^\t X_j, Y_k) \\
&& \hspace{5cm} \times \{\mu_{2p} h_0^2  f( X_j)\}^{-1}
K_{h_{_0}}(X_{ij}) X_{ij}^\t \varepsilon_{b_{_0},i}(Y_k).
\endn
By (C5) and Lemma \ref{basic2}, we have $ n^{-3}\sum_{i,j,k=1}^n
\tilde S_{ijk} =\O\{ (\delta_n + \delta_{ph_0b}^2 +
\delta_n^2/b_0^2  )/h_0\}. $ Thus,
\beginy
n^{-3}\sum_{i,j,k=1}^n  S_{ijk} = n^{-3}\sum_{i,j,k=1}^n \tilde
S_{ijk} + \O\{ r_{ph_{_0}b_{_0}}\delta_{ph_{_0}b_{_0}} h_0^{-1}\}
= \O(\tilde\lambda_n^{(1)}), \label{eq64}
\endy
where
$$
\tilde \lambda_n^{(1)} = \delta_n/h_0 + \delta_{ph_0b}^2/h_0 +
\delta^2_n/(b_0^2 h_0) +
h_0^4+r_{ph_{_0}b_{_0}}\delta_{ph_{_0}b_{_0}} h_0^{-1}.
$$
By (C3) and the strong law of large numbers for U-statistics (cf.
Hoeffding, 1961),  $ \Lambda_n^{(1)} = \int \rho(f(x))
\rho(f_{_Y}(y)) \btd g_{b_0} (B_0^\top x, y) \btdt g_{b_0}
(B_0^\top x, y) \} f(x) $ $f_{_Y}(y) dxdy + o(1)$ almost surely,
which is of full rank asymptotically. Thus its eigenvalues are
greater than a positive constant asymptotically. On the other
hand, the eigenvalues of the lower right principal submatrix in $
C_n $ are of order $ \tilde \lambda_n^{(1)} $. Let
$\lambda_1^{(1)}\ge ... \ge \lambda_p^{(1)} $ be the eigenvalues
of $ \hat \Sigma_{(1)} $ and $ \beta_1^{(1)}, \cdots,
\beta_p^{(1)} $ be the corresponding eigenvectors. By the
interlacing theorem (cf. Ando, 1987), we have $
\min\{\lambda^{(1)}_1, \cdots, \lambda^{(1)}_q\}>c$ and $
\max\{\lambda^{(1)}_{q+1}, \cdots, \lambda^{(1)}_p\} = \O(\tilde
\lambda_n^{(1)})$. By (\ref{opga}) and (\ref{eq64}) we have
\beginy
 \hat \Sigma_{(1)} =  B_0 \Lambda_n^{(1)} B_0^\top +
\O(\delta_B^{(1)}) ,  \label{eq65}
\endy
where $\delta_B^{(1)} = r_{ph_0 b_0} + \delta_{ph_0 b_0} +
\delta_{ph_0 }^2/h_0^2 + \delta_n/h_0 + \delta^2_n/(b_0^2 h_0)$.
Let $ B_{(1)} = (\beta_1^{(1)}, ..., \beta_1^{(q)})$. By Lemma 3.1
of Bai et al (1991), we have
\beginy
 B_{(1)} B_{(1)}^\top - B_0 B_0^\top = \O(\delta_B^{(1)}). \label{rate0}
\endy
Let $ t = 1 $. Consider the $(t+1)$th iteration. Let $
\E_{n,0}^{(t)}(x,y) = \E_{n,0}^{\hat \Sigma_{(t)}}(x,y) $ as
defined in Lemma \ref{krqp}. By the conditions on bandwidths in
(C5), we have $ e_n^{(1)} \stackrel{def}=
\tilde\lambda_n^{(1)}/h_1^2 \to 0 $ and $ \delta_B^{(1)}/h_1 \to 0
$. By Lemma \ref{krqp}, similar to (\ref{opga}), we have from the
algorithm
\beginy
 \hat \Sigma_{(t+1)} = (B_0, \tilde B_0) C^{(t)}_n
(B_0, \tilde B_0)^\top +
 n^{-2} \hspace{-0.1cm}\sum_{j,k=1}^n\{ S^{(t)}_{jk}
 + (S^{(t)}_{jk})^\top\} + \O(\epsilon_{qh_tb_t}\delta_{qh_tb_t}), \label{opgb}
\endy
where $  S^{(t)}_{jk} = \rho_{jk}^{(t)}  B_0 \{\btd g_{b_{_t}}(
B_0^\t X_j, Y_k) + \O(h^2 + \delta_{qh_t}+e^{(t)}_n)\}
\{\E_{n,0}^{(t)}(X_j, Y_k)\}^\top $ and
$$
C^{(t)}_n = \left(\begin{array}{cc} \Lambda_n^{(t)} & \O(\epsilon_{qhb}) \\
\O(\epsilon_{qhb}) & \O(\epsilon^2_{qhb})
\end{array}\right),
$$
where $ \Lambda_n^{(t)} = n^{-2}\sum_{j,k=1}^n \rho_{jk}^{(t)}
\btd g_{b_{_t}}( B_0^\t X_j, Y_k)  \btdt g_{b_{_t}}( B_0^\t X_j,
Y_k) + \O\{h_t^2 + \delta_{qh_t} +
  e_n^{(t)}\} $.  Note that $ B_{(t)}^\top \E_{n,0}^{(t)}(X_j, Y_k) = 0
$, $ \E_{n,0}^{(t)}(X_j, Y_k) = \O(\delta_{qh_tb_t})$ and  $
B_0^\top \E_{n,0}^{(t)}(X_j, Y_k) = \O(
\delta_{qh_tb_t}\delta_B^{(t)})$. It follows that
\beginy
\n &&  n^{-2} \hspace{-0.1cm}\sum_{j,k=1}^n\{ S^{(t)}_{jk}
 + (S^{(t)}_{jk})^\top\}\\
\n &&= (B_0, \tilde B_0)\Big[ (B_0, \tilde B_0)^\top n^{-3}
\hspace{-0.1cm}\sum_{j,k=1}^n\{ S^{(t)}_{jk}
 + (S^{(t)}_{jk})^\top\} (B_0, \tilde B_0) \Big]  (B_0, \tilde
B_0)^\top\\
&&= (B_0, \tilde B_0)\hspace{-0.1cm}\left(\hspace{-0.2cm}\begin{array}{cc} 0 & C^{(t)}_{12,n} \\
(C^{(t)}_{12,n})^\top & 0
\end{array}\hspace{-0.2cm}\right) (B_0, \tilde B_0)^\top+\O(
\delta_{qh_tb_t}\delta^{(t)}_B), \qquad \label{eq68}
\endy
where $ C^{(t)}_{12,n} = n^{-2} \sum_{j,k=1}^n
\rho_{jk}^{(t)}\{\btd g_{b_t}(B_0^\top X_j, Y_k)+ \O(h_t^2 +
\delta_{qh_t}+e^{(t)}_n)\}\{\E_{n,0}^{(t)}(X_j, Y_k)\}^\top $
$\tilde B_0 $. Similar to $ \rho^{(0)}_{xy} $, we have $
\rho_{jk}^{(t)} = \tilde \rho_{jk}^{(t)} + \O(r_{qh_tb_t}) $ where
$ \tilde \rho_{jk}^{(t)} = \rho(f_{B_0}(X_j))\{\rho(f_{_Y}(Y_k))$
$+\rho'(f_{_Y}(Y_k))f_{_Y}''(Y_k)b_t^2/2\} $. By (C5) and Lemma
\ref{basic2}, we have
\beginy
\n C^{(t)}_{12,n} \eqc n^{-2} \hspace{-0.1cm}\sum_{j,k=1}^n \tilde
\rho_{jk}^{(t)}\btd g_{b_t}(B_0^\top X_j,
Y_k)\{\E_{n,0}^{(t)}(X_j, Y_k)\}^\top \tilde B_0 + \O(r_{qh_tb_t}
\delta_{qh_tb_t} + e^{(t)}_n
\delta_{qh_tb_t})\\
\eqc \O(\delta_n + \delta_{qh_t b_t}^2 + \delta_n^2 b_t^{-2} +
r_{qh_tb_t} \delta_{qh_t b_t}+ e^{(t)}_n \delta_{qh_tb_t}).
\label{cc12}
\endy
By the strong law of large numbers for U-statistics, it follows $
\Lambda_n^{(t)} = M_0 + o(1) $ almost surely, where $ M_0 $ is
defined in (C3). Let $ \lambda_1^{(t+1)} \ge ... \ge
\lambda_p^{(t+1)} $ be the eigenvalues of $  \hat \Sigma_{(t+1)} $
and $ B_{(t+1)} $ the first $ q$ eigenvectors. By the same
arguments as for $\tilde \lambda_n^{(1)} $, it follows from
(\ref{opgb}), (\ref{eq68}) and (\ref{cc12}) that $ \min\{
\lambda_1^{(t+1)}, ..., \lambda_q^{(t+1)} \}
> c$ and $ \max\{ \lambda_{q+1}^{(t+1)}, ..., \lambda_p^{(t+1)} \}
= \O\{\tilde \lambda_n^{(t+1)} \} $, where $ \tilde
\lambda_n^{(t+1)} = \epsilon_{qh_tb_t} \delta_{qh_tb_t} +
\epsilon_{qh_tb_t}^2+ \delta_{qh_tb_t}\delta_B^{(t)}$. Considering
$ e_n^{(t+1)} h_{t+1}  \stackrel{def}= \tilde
\lambda_n^{(t+1)}/h_{t+1}$, there exists a constant $c_1$, which
does not depend on $ t$, such that
\beginy
 e_n^{(t+1)} h_{t+1} \le c_1\{ \chi_{0,n}^{(t)} +
    \chi^{(t)}_{1,n} e_n^{(t)} h_t +
   \chi_{2,n}^{(t)}\delta_B^{(t)}\},
   \label{rate1}
\endy
where $ \chi_{0,n}^{(t)} = (h_t^4+ h_t^2 \delta_{qh_tb_t} +
\delta_{qh_tb_t}\delta_{qh_t})\delta_{qh_tb_t}/h_{t+1}$, $
\chi_{1,n}^{(t)} = (h_t^2 + \delta_{qh_tb_t}) \delta_{qh_tb_t}
/(h_t h_{t+1})$ and $ \chi_{2,n}^{(t)} =
\delta_{qh_tb_t}/h_{t+1}$. By (\ref{opgb}) and (\ref{eq68}), we
write
\beginy
\hat \Sigma_{(t+1)} = B_0 \Lambda^{(t)}_n B_0 + B_0 \tilde
C^{(t)}_{12,n} \tilde B_0^\top + \tilde B_0 (\tilde C^{(t)}_{12,n}
B_0)^\top+ \O\{\epsilon_{qh_tb_t} +
\delta_{qh_tb_t}\delta_B^{(t)}\}, \label{opgb1}
\endy
where $ \tilde C^{(t)}_{12,n} $ is the first term on the right
hand side of the first equation in (\ref{cc12}). By the same
arguments as for (\ref{rate0}), we have
 $
   B_{(t+1)} B_{(t+1)}^\top - B_0 B_0^\top  = \O\{\delta_{qh_tb_t}( \delta_{qh_tb_t} +
r_{qh_tb_t})+ (h_t^2 + r_{qh_tb_t})e_n^{(t)} +
(h+\delta_{qhb}/h)\delta_B^{(t)}  + \delta_n +\delta_n^2/b_t^2 +
h_t^4\}
 $. That is
\beginy
\delta_B^{(t+1)} \le c_2\{ \chi_{3,n}^{(t)}+ \chi_{4,n}^{(t)}
e_n^{(t)} h_t + \chi_{5,n}^{(t)}\delta_B^{(t)} \} \label{bB}
\endy
for a constant $ c_2 $ independent of $ t$, where $
\chi_{3,n}^{(t)} = \delta_{qh_tb_t}( \delta_{qh_tb_t}+r_{qh_tb_t})
 + h_t^4 + \delta_n^2/b_t^2+ \delta_n$, $ \chi_{4,n}^{(t)} = (h_t^2 +
r_{qh_tb_t})/h_t $ and $ \chi_{5,n}^{(t)} =
h_t+\delta_{qh_tb_t}/h_t$. Note that $ h_t $ and $ b_t $
decreasing with $ t $, by (C5) we have $\delta_{qh_t b_t}/h_{t+1}
\le \delta_{q\hbar\bbar}/\hbar \to 0 $. It follows that $
e_n^{(t+1)} = \lambda_n^{(t+1)}/h_{t+1}^2 \to 0 $, $
\delta_B^{(t+1)} = O(r_{q h_t b_t}) $ and $
\delta_B^{(t+1)}/h_{t+1} \to 0$.  Recursing (\ref{rate1}) and
(\ref{bB}), it follows  that
\beginn
 \delta_B^{(\infty)}  = \O\{\chi_{3,n}^{(\infty)} + \chi_{4,n}^{(\infty)} \chi_{0,n}^{(\infty)}\}
  = \O\{ \hbar ^4 + \delta_{q\hbar  \bbar }( \delta_{q\hbar  \bbar } +
 \hbar ^2  + \bbar ^4) + \delta_n^2/\bbar ^2 + \delta_n\}
\endn
and $ e_n^{(\infty)} = \O(\delta_{q\hbar\bbar}) $. This is the
first part of Theorem \ref{mainOPG}.  By (\ref{opgb1}) and the
equations above, write
\beginn
\hat \Sigma_{(\infty)} &=&  \{B_0 + \eta_n \} \Lambda^{(\infty)}_n
\{B_0 + \eta_n \}^\top +  \O\{\hbar ^4 + \delta_{q\hbar  \bbar }(
\delta_{q\hbar  \bbar } +
  \bbar ^4)+ \delta_n^2/\bbar ^2\},
\endn
where $ \eta_n = \tilde C^{(\infty)}_{12,n}
(\Lambda_n^{(\infty)})^{-1} = \O\{\hbar ^4 + \delta_{q\hbar  \bbar
}( \delta_{q\hbar  \bbar }  + \bbar ^4) + \delta_n^2/\bbar ^2 +
\delta_n\} $.
Note that $ B_{(\infty)}^\top \bar w_{_{B_{(\infty)}}}^+(x) = 0 $
and thus $ B_{(\infty)}^\top \eta_n = 0 $. We have $ \tilde
\Lambda_n \stackrel{def}= (B_0+\eta_n)^\top (B_0+\eta_n) = I_q +
\O( \delta_n^2). $
Let $ \tilde \eta_n = \{B_0 + \eta_n \} \tilde \Lambda_n^{-1/2} $.
It follows that
\beginn
\hat \Sigma_{(\infty)} &=&   \tilde \eta_n \Lambda^{(\infty)}_n
\tilde \eta_n^\top +  \O\{\hbar ^4 + \delta_{q\hbar  \bbar }(
\delta_{q\hbar \bbar } +
  \bbar ^4)+ \delta_n^2/\bbar ^2\}.
\endn
Let $ \hat B_{dOPG} $ be the first $q$ eigenvectors of $\hat
\Sigma_{(\infty)} $. By Lemma 3.1 of Bai et al (1991), we have
\beginy
\hat B_{dOPG} \hat B_{dOPG}^\top  - B_0 B_0^\top = B_0 \eta_n^\top
+ \eta_n B_0^\top + O\{ \hbar ^4 + \delta_{q\hbar \bbar }(
\delta_{q\hbar \bbar } + \bbar ^4)+ \delta_n^2/\bbar ^2\}. \ \
\label{normal1}
\endy
By Lemma \ref{basic2} and (C5), we have
\beginy
\n \eta_n &=& n^{-2}\sum_{j,k=1}^n
 \rho(f_{B_0}(X_j))\rho(f_{_Y}(Y_k)) \E_{n,0}^{(\infty)}(X_j, Y_k) \btdt g_b(B_0^\t X_j,
Y_k) (\Lambda_n^{(\infty)})^{-1} \\
\n && \hspace{10cm}+ \O\{ r_{q\hbar  \bbar }  \delta_{q\hbar  \bbar }\} \quad \ \\
\n &=& n^{-1}
 \sum_{i=1}^n\rho(f_{B_0}( X_i))\rho(f_{_Y}(Y_i)) \bar w^+_{_{B_0}}(X_i)
 \nu_{_{B_0}}(X_i)\zeta^\t_i (\Lambda_n^{(\infty)})^{-1}
 + \O\{ r_{q\hbar  \bbar }  \delta_{q\hbar  \bbar } \},
\endy
where  $ \ \zeta_i  \ = \ \btd g_\bbar (B_0^\t X_i, Y_i)
f_{_Y}(Y_i) - E\{\btd g_\bbar (B_0^\t X_i, Y_i )
f_{_Y}(Y_i)|B_0^\t X_i\} $. \ Let $\tilde \zeta_i = $ $ \btd
f(Y_i|B_0^\t X_i) f_{_Y}(Y_i) $ $- E\{\btd f(Y_i|B_0^\t X_i)
f_{_Y}(Y_i)|B_0^\t X_i\} $. As $ b \to 0 $, we have $
\Lambda_n^{(\infty)} \to M_0 $ almost surely, where $ M_0 $ is
defined in (C3). By calculating the mean and covariance matrix, we
have
\beginn
    n^{-1}
 \sum_{i=1}^n\rho(f_{B_0}( X_i))\rho(f_{_Y}(Y_i)) \bar w^+_{_{B_0}}(X_i)
 \nu_{_{B_0}}(X_i)(\tilde
 \zeta^\top_i - \zeta^\t_i)   =
 o_P(n^{-1/2}).
\endn
It follows from the two equations above and the conditions in the
Theorem for the bandwidths
\beginy
 \eta_n =  n^{-1}
 \sum_{i=1}^n\rho(f_{B_0}( X_i))\rho(f_{_Y}(Y_i)) \bar w^+_{_{B_0}}(X_i)
 \nu_{_{B_0}}(X_i)\tilde \zeta^\t_i M_0^{-1}
 + o_P(n^{-1/2}). \quad \label{norma2}
\endy
After vectorizing $ \eta_n $, the second part of Theorem
\ref{mainOPG} follows from (\ref{normal1}), (\ref{norma2}) and the
central limit theorem. \hspace{\fill}$ \Box $

{\bf Proof of Theorem \ref{mainMAVE}} Consider the initial
estimator $  B_{(1)} $ in (\ref{rate0}). Let $ \tilde Q =
B_{(1)}^\top B_0$. For simplicity, we assume $ \tilde Q = I_q $;
otherwise, we may use basis $ B_0 \tilde Q $ and consider the
expansion in Lemmas \ref{krq}, \ref{denominator} and
\ref{numerator} at $ (B_0\tilde Q)^\top x$. Let $ \tilde
\delta_B^{(t)} $ be the consistency rate of the estimator in the
$t'$th iteration. Write $
 \ell(B_0) = (I_q\otimes B_0)\ell(I_q).
$  By the definition of $ D_B $ in Lemma \ref{denominator}, it
follows
\beginy
(I_q\otimes B)^\top D_B = 0, \quad I_q\otimes B = I_q\otimes B_0 +
O(\delta_B), \quad (I_q\otimes B_0)^\top \Phi_n(B_0) = 0. \ \
\label{kfngkedf}
\endy
By the definition of the Moore-Penrose inverse we have  $ D_B^+
D_B = I_q\otimes (\tilde B \tilde B^\top) ,$ where $ (B, \tilde B)
$ is a $ p\times p $ orthogonal matrix. By Lemmas
\ref{denominator}, \ref{numerator} and (\ref{relation0}), for
every $ B_{(t)} $ in $ {{\cal B}} = \{ B: |B-B_0| \le \tilde
\delta_B^{(t)}$\}, if $ \tilde \delta_B^{(t)}/h_{t} \to 0 $ we
have
\beginy
\n {\bf b}^{(t+1)}  &=& (I_q\otimes B_0)\{\ell(I_q)+O(c^{(t)}_n)\}
+ \frac12 \Psi_{(t)} \{ \ell(B_{(t)}) - \ell( B_0)\}+ \frac12
D_{(t)}^{+} \Phi_n(B_0)
  \\
&&  + \O\{ \Delta_t + (h_t+ \delta_{qh_tb_t}/h_t) \tilde
\delta_B^{(t)} \}, \quad \label{fff}
\endy
where $\Delta_t = h_t^4 + (h_t^2 + b_t^4 +
\delta_{qh_tb_t})\delta_{qh_{t}b_{t}} + \delta_n^2/b_t^2$, $
c^{(t)}_n = \{\Delta_t +  (\delta_{qh_tb_t}/h_t+h_t)\tilde
\delta_B^{(t)}\}/h_t^2$, $D_{(t)} =  D_{B_{(t)}} $ and $
\Psi_{(t)} = I_q\otimes (\tilde B_{(t)} \tilde B_{(t)}^\top) =
\Psi + \tilde \delta_B^{(t)} $, where $ \Psi = I_q\otimes ( \tilde
B_0 \tilde B_0^\top)$ is a projection matrix and $ (B_0, \tilde
B_0) $ is
 a $ p\times p $ orthogonal matrix.  We have
\beginn
  {{\cal M}}( {{\bf
b}}^{(t+1)})  &=& B_0 \Lambda_n^{(t)} + \frac12 {{\cal
 M}}( \Psi \{\ell(B_{(t)})-\ell( B_0)\})
 +\frac12 {{\cal
 M}}(D_{(t)}^{+} \Phi_n(B_0) )\\
 && \hspace{6cm} +
\O\{ \Delta_t + (h_t+ \delta_{qh_tb_t}/h_t) \tilde
\delta_B^{(t)}\},
\endn
where $ \Lambda_n^{(t)} = I_q + O(c^{(t)}_n) $ and $ {{\cal M}}(.)
$ is defined in section 2.2. Note that
$$ \tilde \Lambda_n^{(t+1)}
\stackrel{def}= \{{{\cal M}}( {{\bf b}}^{(t+1)})\}^\top {{\cal
M}}( {{\bf b}}^{(t+1)}) = (\Lambda^{(t)}_n)^2 + \O\{
\delta_B^{(t)} + \tilde \delta_n + \Delta_t + (h_t+
\delta_{qh_tb_t}/h_t) \tilde \delta_B^{(t)}\},
$$
where $ \tilde
\delta_n = \delta_n + \delta_{qh_tb_t}^2/h_t $. If  $
 c_n^{(t)} = o(1) $ almost surely, then by Step 3
\beginy
\n  B_{(t+1)}\hspace{-0.2cm} &=&\hspace{-0.2cm}  B_0 +  \frac12
{{\cal
 M}}( \Psi \{\ell(B_{(t)})-\ell( B_0)\})+\frac12 {{\cal
 M}}(D_{(t)}^{+} \Phi_n(B_0) ) \\
\n && \hspace{7cm} +\O\{  \Delta_t + (h_t+ \delta_{qh_tb_t}/h_t)
\tilde
\delta_B^{(t)}\} \ \\
\hspace{-0.2cm}&=&\hspace{-0.2cm}  B_0 + \frac12 {{\cal
 M}}(
  \Psi \{ \ell( B_{(t)})-\ell( B_0)\}) +\O\{ \tilde \delta_n + \Delta_t + (h_t+ \delta_{qh_tb_t}/h_t) \tilde
\delta_B^{(t)}\}. \label{fd}
\endy
By (C5) and (\ref{rate0}), we have $\delta_{qh_t b_t }/h_t^2 \le
\delta_{q\hbar \bbar }/\hbar^2 \to 0$, $ \delta_B^{(1)}/h_1 \to 0$
and $ c_n^{(1)} \to 0 $ almost surely. Thus (\ref{fd}) holds for $
t = 1$. By assumption (C5), it follows that $ \tilde
\delta_B^{(2)}/h_{2} =o(1) $ and $ c_n^{(2)} = o(1) $ almost
surely. Thus (\ref{fd}) holds for $ t = 2$.   Recurring the
formula, we have
 $$
\tilde \delta_B^{(\infty)} = \O(\Delta_\infty +\tilde \delta_n ) =
\O\{ \hbar ^4 + (\hbar ^2+ \bbar ^4 +\delta_{q\hbar  \bbar }
)\delta_{q\hbar \bbar }+\tilde \delta_n \}.
$$
A more detailed deduction was given in Xia, Tong and Li  (2002).
Therefore, the first part of Theorem \ref{mainMAVE} follows
immediately.  By the first equation of (\ref{fd}) with $ t =
\infty$ and Lemma \ref{numerator}, we have
\beginn
\n B_{(\infty)} - B_0 &=&  \frac12 {{\cal M}}(\Psi \{ \ell(
B_{(\infty)}) - \ell( B_0)\}) + \frac12 {{\cal M}}(D_{(\infty)}^+
\Phi_n(B_0)) \qquad \qquad
\\
&& \hspace{3.5cm}
 + O_P\{\hbar ^4 + (\hbar ^2+ \bbar ^4 +\delta_{q\hbar  \bbar } )\delta_{q\hbar  \bbar }
 \}.
\endn
Multiplying both sides by $ B^\top_0$, by (\ref{kfngkedf}) we have
\beginn
B_0^\top B_{(\infty)} - I&=&  O_P\{\hbar ^4 + (\hbar ^2+ \bbar ^4
+\delta_{q\hbar  \bbar } )\delta_{q\hbar  \bbar }
 \}.
\endn
It follows that
\beginn
\n B_{(\infty)}B_{(\infty)}^\top B_0 - B_0 &=&  \frac12 {{\cal
M}}(\Psi \{ \ell( B_{(\infty)}) - \ell( B_0)\}) + \frac12 {{\cal
M}}(D_{(\infty)}^+ \Phi_n(B_0)) \qquad \qquad
\\
&& \hspace{3.5cm}
 + O_P\{\hbar ^4 + (\hbar ^2+ \bbar ^4 +\delta_{q\hbar  \bbar } )\delta_{q\hbar  \bbar }
 \}.
\endn
 Note that $ \Psi D_{(\infty)}^+ = D_{(\infty)}^+
+ O_P(\tilde \delta_B^{(\infty)})$. We have
\beginn
\ell( B_{(\infty)} B_{(\infty)}^\top B_0)  - \ell( B_0)  =
D_{(\infty)}^+ \Phi_n(B_0) + O_P\{\hbar ^4 + (\hbar ^2+ \bbar ^4
+\delta_{q\hbar  \bbar } )\delta_{q\hbar  \bbar } \}.
\endn
This is the second part of  Theorem \ref{mainMAVE}. \hspace{\fill}
$ \Box $


\subsection{Auxiliaries}

\begin{Lemma}\label{basic1} \rm
Suppose $m_n(\chi, Z), n = 1, 2, \cdots, $ are measurable
functions of $Z $ with index $\chi \in {\Bbb{R}}^d $, where $d $
is an integer, such that (I) $| m_n(\chi, Z) | \le  M(Z) $ with $E
(M^r(Z)) < \infty $ for some $ r > 2$; (II)  $\sup_{\chi}
E|m_n(\chi, Z)|^2 < a_n $; and (III) $|m_n(\chi,Z) -
m_n(\chi^{\prime}, Z) | \le |\chi-\chi^{\prime}|^{\alpha_1}
n^{\alpha_2} G(Z)$ with some $\alpha_1, \alpha_2 >0 $ and $E|G(Z)|
 < \infty $. Suppose $\{Z_i, i = 1, \cdots, n\} $ is a random sample
from $ Z$. If $a_n = c n^{-\delta}$ with $ 0 \le \delta < 1-2/r$
and $ c > 0 $, then for any positive $\alpha_0 $ we have
\begin{eqnarray*}
\sup_{|\chi| \le n^{\alpha_0} }\Big| n^{-1} \sum_{i=1}^n\{
m_n(\chi, Z_i) - E m_n(\chi, Z_i)\}\Big| = O\{ ( a_n \log n
/n)^{1/2}\}
\end{eqnarray*}
almost surely.
\end{Lemma}

{\bf Proof of Lemma \ref{basic1}}  The ``continuity argument''
approach is used here. See, e.g. Mack and  Silverman (1982) and
H\"ardle et al  (1993). Note that $ \D _n \stackrel{def}{= }\{
|\chi| \le n^{\alpha_0} \} $ is bounded and its Borel measure is
less than $c_1 n^{\alpha_0 d} $ for some constant $ c_1$. There
are $ n^{\alpha_4} $ $( \alpha_4
>   \alpha_0 d + (1+\alpha_2)d/\alpha_1 $) balls $B_{n_k} $ centered at $ \chi_{n_k} $,
$1\le k \le n^{\alpha_4} $, with diameter less than $c_2
n^{-(1+\alpha_2)/\alpha_1} $, such that $\D _n \subset \cup_{ 1
\le k \le n^{\alpha_4} } B_{n_k} $. It follows that
\begin{eqnarray}
\n && \sup_{\chi \in \D _n }\Big|\frac{1}{n}\sum_{i=1}^n\{
m_n(\chi, Z_i) - Em_n(\chi, Z_i) \}\Big|\\
\n &&\le \max_{ 1 \le k \le n^{\alpha_4} }
\Big|\frac{1}{n}\sum_{i=1}^n\{ m_n(\chi_{n_k}, Z_i) - E
m_n(\chi_{n_k}, Z_i) \}\Big| \\
\n && \qquad   + \max_{ 1 \le k \le n^{\alpha_4} } \sup_{ \chi \in
B_{n_k} } \Big| \frac{1}{n}\sum_{i=1}^n \Big[ \{
m_n(\chi, Z_i) - m_n(\chi_{n_k}, Z_i) \} \\
&& \hspace{1cm}
 - E\{ m_n(\chi, Z_i) - m_n(\chi_{n_k}, Z_i)\}\Big]
\Big|
\nonumber \\
&& \stackrel{def}{=} \max_{ 1 \le k \le n^{\alpha_4} } |R_{n,k,1}|
+ \max_{ 1 \le k \le n^{\alpha_4} } \sup_{ \chi \in B_{n_k} }
|R_{n,k,2}|. \label{ffaacc0}
\end{eqnarray}
By condition (III) and the definition of $ B_{n_k} $, we have
\begin{eqnarray*}
\ \ \max_{ 1\le k \le n^{\alpha_4} } \sup_{ \chi \in B_{n_k}} |
m_n(\chi, Z_i) -m_n(\chi_{n_k}, Z_i) | &\le& \max_{ 1\le k \le
n^{\alpha_4}} \sup_{\chi \in B_{n_k}} n^{\alpha_2}|\chi -
\chi_{n_k} |^{\alpha_1}G(Z_i) \\
&\le& c_3 n^{-1} G(Z_i).
\end{eqnarray*}
By the strong law of large numbers, we have
\begin{eqnarray}
\max_{ 1 \le k \le n^{\alpha_4} } \sup_{ \chi\in B_{n_k} }
|R_{n,k,2}| \le c_4 n^{-2} \sum_{i=1}^n \{G(Z_i) + E G(Z_i)\} =
O(n^{-1}) \label{ffaacc1}
\end{eqnarray}
almost surely.  Let $T_n = (n a_n/\log n)^{1/2} $,
$m_{n}^o(\chi_{n_k}, Z_i) = m_n(\chi_{n_k}, Z_i) I\{ |M(Z_i)|$
$\ge T_n\} $ and $m_{n}^I(\chi_{n_k}, Z_i) = m_n(\chi_{n_k}, Z_i)
- m_{n}^o(\chi_{n_k}, Z_i)$. Write
\begin{eqnarray}
R_{n,k,1} =\frac{1}{n}\sum_{i=1}^n\Big[ m_{n}^o(\chi_{n_k}, Z_i) -
E\{m_{n}^o(\chi_{n_k}, Z_i) \}\Big] +
\frac{1}{n}\sum_{i=1}^n\xi_{n_k, i}, \label{rnk1}
\end{eqnarray}
where $\xi_{n_k, i} = m_n^I(\chi_{n_k}, Z_i) -
E\{m_n^I(\chi_{n_k}, Z_i)\} $. By the truncation, it follows that
$$
E | m_{n}^o(\chi_{n_k}, Z_i)| \le T_n ^{-r+1} E |M(Z_i)|^r.
$$
If $a_n = c n^{-\delta} $ with $0 \le \delta < 1 - 2/r$, we have
\begin{eqnarray}
n^{-1} |\sum_{i=1}^{n} E m_{n}^o(\chi_{n_k}, Z_i) |  \le  E
|M(Z_1)|^r T_n^{-r+1} =o (
\{ a_n \log(n)/n\}^{1/2}). \label{df1} 
\end{eqnarray}
Again by the truncation, we have
$$ 
\sum_{i=1}^n
|m_{n}^o(\chi_{n_k}, Z_i)| \le \sum_{i=1}^n |M(Z_i)|I( |M(Z_i)|
\ge T_n) \le T_n^{-r+1} \sum_{i=1}^n |M(Z_i)|^r I( |M(Z_i)| \ge
T_n).
$$ 
For fixed $T$, by the strong law of large
numbers, we have
\begin{eqnarray*}
n^{-1} \sum_{i=1}^n |M(Z_i)|^r I( |M(Z_i)| \ge
T ) \to E\{ |M(Z_1)|^r I( |M(Z_1)| \ge T ) \}
\end{eqnarray*}
almost surely. The right hand side above is dominated by $E\{
|M(Z_i)|^r \} $ and $\to 0 $ as $T \to \infty $. Note that $T_n $
increase to $\infty$ with $ n $. For large $n $ such that $T_n > T
$, we have
\begin{eqnarray*}
 \ C_n \stackrel{def}= n^{-1} \sum_{i=1}^n |M(Z_i)|^r I( |M(Z_i)| \ge
T_n) \le n^{-1} \sum_{i=1}^n |M(Z_i)|^r I( |M(Z_i)| \ge T ) \to 0
\end{eqnarray*}
almost surely as $T \to \infty $. It follows
\beginy
 \max_{ 1\le k
\le n^{\alpha_4} } n^{-1}|\sum_{i=1}^n m_{n}^o(\chi_{n_k}, Z_i)|
\le C_n T_n^{-r+1}
 = o\{ (a_n \log n/n)^{1/2} \}
\label{df2}
\endy
almost surely. By condition (II), we have
\beginy
\n \max_{1\le k\le n^{\alpha_4}} \mbox{Var}(\sum_{i=1}^n
\xi_{n_k,i})
 &\le & n  \max_{1\le k\le n^{\alpha_4}} E\{ m_{n}^I(\chi_{n_k},
Z_1)\}^2 \\
&\le& n \max_{1\le k\le n^{\alpha_4}} E \{m_{n}(\chi_{n_k},
Z_1)\}^2 \le c_5 n a_n \stackrel{def}{= }N_1. \qquad \label{var1}
\endy
By the  condition on $ a_n $ and the definition of $ \xi_{n_k,i}
$, we have
\begin{eqnarray}
\max_{1\le k\le n^{\alpha}} |\xi_{n_k, i}|\le c_6 T_n = c_6 (n a_n
/\log n)^{1/2} \stackrel{def}{= }N_2.\ \ \label{bound1}
\end{eqnarray}
Let $N_3 = c_7 (n a_n\log n)^{1/2}$ with $ c_7^2 > 2 (\alpha_4 +2)
(c_5+c_6c_7) $. By the Bernstein's inequality (cf. DE LA Pe\~na,
1999), we have from (\ref{var1}) and (\ref{bound1}) that
\begin{eqnarray}
\n P(|\sum_{i=1}^n \xi_{n_k,i} |> N_3) &\le& 2 \exp\left(
\frac{-N_3^2}{2( N_1 + N_2 N_3) } \right)\\
\n  &\le& 2\exp\{ - c_7^2 \log n/(2c_5+ 2c_6c_7)\} \\
\n &\le& c_8 n^{-\alpha_4-2}.  \label{lemma407}
\end{eqnarray}
It follows that
\begin{eqnarray}
 \sum_{n=1}^{\infty}\Pr( \max_{ 1 \le k \le n^{\alpha_4} } |
\sum_{i=1}^n \xi_{n_k,i} | \ge N_3 ) \le \sum_{n=1}^{\infty}
n^{\alpha_4} \max_{ 1 \le k \le n^{\alpha_4} } \Pr( | \sum_{i=1}^n
\xi_{n_k,i} | \ge N_3 )  < \infty. \ \  \label{eqR2a}
\end{eqnarray}
By the Borel-Cantelli lemma (cf. Chow and Teicher, 1978, p.60), we
have
\begin{eqnarray}
\max_{ 1 \le k \le n^{\alpha_4} } |\sum_{i=1}^n \xi_{n_k,i}| =
O(N_3) \quad \label{lemma413}
\end{eqnarray}
almost surely. Combining (\ref{rnk1}), (\ref{df1}), (\ref{df2})
and (\ref{lemma413}), we have
\begin{eqnarray}
\max_{ 1 \le k \le n^{\alpha_4} } |R_{n,k,1}| = O\{ (a_n \log(
n)/n)^{1/2} \}  \label{eqR2}
\end{eqnarray}
almost surely.  Lemma \ref{basic1} follows from (\ref{ffaacc0}), (\ref{ffaacc1}) and (%
\ref{eqR2}). \hspace{\fill}  $ \Box $

For any function $ G(X_i, Y_i, X_j,$ $ Y_j, X_k, Y_k) $ (or $
G(X_j,$ $ Y_j, X_k, Y_k) $), we introduce a projection operator $
E_k $ as follows.
\beginn
E_k G(X_i, Y_i, X_j, Y_j, X_k, Y_k) = E\{ G(X_i, Y_i, X_j, Y_j,
X_k, Y_k) |X_i, Y_i, X_j, Y_j\}. \label{operator}
\endn

\begin{Lemma}\label{basic2} \rm Let $ {{\cal A}} = \{A: A^\t A =
I_{\kappa}\} $ with $ 1\le \kappa \le p $. Suppose $ g_0(y),
g_1(x), g_2(x) $ are bounded continuous functions. If conditions
(C2) and (C4) hold with $ B $ replaced by $ A$ for all $ A \in
{{\cal A}}$, then
\beginn
 n^{-3} \sum_{i,j,k=1}^n K_h (A^\t X_{ij})
  g_1(X_i) g_2(X_j) g_0(Y_k) \btd g_b(B_0^\top X_j, Y_k)
  \varepsilon_{b,i} (Y_k) \hspace{3cm}\\
 = n^{-1}\sum_{i=1}^n E_j E_k \{ K_h( A^\top X_{ij}) \btd
g_b(B_0^\top X_j, Y_k)\varepsilon_{b,i}(Y_k)\} +
\O(\varsigma_{\kappa hb}|A \in {{\cal A}}),
\endn
where $ \varsigma_{\kappa hb} = \delta_n^3  h^{-\kappa}b^{-2} +
\delta_{\kappa h b}^2 +  \delta_n^2 b^{-2}$ and the first term on
the right hand side is $ \O(\delta_n)$.
\end{Lemma}

\noindent {\bf Proof of Lemma \ref{basic2}} For easy of
exposition, we consider $ g_k \equiv 1, k = 0, 1, 2 $ only. Let $
\Delta_n(A) $ be the left hand side of the equation in the lemma.
Let $\varphi_{_K}(s)  =  (2\pi)^{-\kappa}$ $\int \exp(\i s^\top u)
K(u)du $ and $ \varphi_{_H}(t) =(2\pi)^{-1} \int \exp(\i t v)
H(v)dv $ be the Fourier transformations, where $\i $ is the
imaginary unit. It follows from the inverse Fourier transformation
that $ g_b(u, y) = b^{-1} \int \varphi_{_H}(t') e^{-\i t' y/b} E\{
e^{\i t' Y/b}| $ $B_0^\top X = u\} dt' $. Thus
\beginy
 \btd g_b(B_0^\top X_{j}, Y_k) = b^{-1} \int
\varphi_{_H}(t') \btd \tilde g_b(B_0^\top X_j) e^{-\i t'Y_k/b}
dt', \label{cccai}
\endy
 where $ \btd \tilde g_b(u) =
\partial E(e^{\i t'Y/b}|B_0^\top X=u)/\partial u $. We have
\beginy
\n \Delta_n(A ) \eqc \frac{1}{n^3  b}\int \varphi_{_H}(t')
\sum_{i,j,k=1}^n\{ K_h( A ^\top X_{ij})  \btd \tilde g_b(B_0^\top
X_j) - E_j [K_h( A ^\top
X_{ij}) \\
\n && \hspace{1.5cm} \times \btd \tilde g_b(B_0^\top X_j)]\}
\{\varepsilon_{b,i}(Y_k) e^{-\i t'Y_k/b} - E_k[\varepsilon_{b,i}(Y_k) e^{-\i t'Y_k/b}]\}dt'\\
\n &&\hspace{-0.2cm} + \frac{1}{n^2 b} \int \varphi_{_H}(t')
\sum_{i,k=1}^n E_j [K_h( A ^\top
X_{ij})  \btd \tilde g_b(B_0^\top X_j)]\\
\n && \hspace{3.8cm} \times \{\varepsilon_{b,i}(Y_k) e^{-\i
t'Y_k/b} -
E_k[\varepsilon_{b,i}(Y_k) e^{-\i t'Y_k/b}]\} dt'\\
\n &&\hspace{-0.2cm}  + \frac{1}{n^2 b } \int \varphi_{_H}(t')
\sum_{i,j=1}^n  E_k[\varepsilon_{b,i}(Y_k) e^{-\i t'Y_k/b}] \{
K_h(
A ^\top X_{ij})  \btd \tilde g_b(B_0^\top X_j) \\
\n && \hspace{5.4cm} - E_j [K_h( A ^\top
X_{ij})  \btd \tilde g_b(B_0^\top X_j)]\}dt'\\
\n &&\hspace{-0.2cm}  + \frac{1}{n b } \int \varphi_{_H}(t')
\sum_{i=1}^n  E_j
[K_h( A ^\top X_{ij})  \btd \tilde g_b(B_0^\top X_j)] E_k[\varepsilon_{b,i}(Y_k) e^{-\i t'Y_k/b}]dt'\\
 &\stackrel{def}=& \Delta_{n,1}(A )
+ \Delta_{n,2}(A ) +\Delta_{n,3}(A ) + \Delta_{n,4}(A ).
\label{adf07}
\endy
By the inverse Fourier transformation, it follows that  $ K_h(A
^\top X_{ij}) = h^{-\kappa} \int \varphi_{_K}(s) $ $e^{-\i s^\top
A ^\top X_{ij}/h}ds  $ and $ H_b(Y_i-Y_k) = b^{-1} \int
\varphi_{_H}(t) e^{-\i t (Y_i-Y_k)/b}dt  $ .
 Thus
\beginn
\Delta_{n,1}(A ) &=&  \frac{1}{n^3h^{\kappa}b^2} \int
\prod_{\ell=1}^3 \sum_{i=1}^n
           m_{\ell,n}(A, s, t, t', X_i, Y_i) \varphi_{_K}(s) \varphi_{_H}(t) \varphi_{_H}(t') ds dt
              dt',
\endn
where
$$ m_{1,n}(A, s, t, t', X_i, Y_i) =
            e^{ - \i s^\top A ^\top X_i /h}\btd \tilde g_b(B_0^\top X_i) -
            E [e^{ - \i s^\top A ^\top X_i /h}\btd \tilde g_b(B_0^\top
            X_i)],
             $$
$$             m_{2,n}(A, s, t, t', X_i, Y_i) =  e^{ \i (t-t')Y_i/b} - E[e^{
\i( t-
             t')Y_i/b}]
$$
and
$$
m_{3,n}(A, s, t, t', X_i, Y_i) =  e^{- \i t
            Y_i/b} - E(e^{- \i t
            Y_i/b}|X_i) .
$$

By (C2), we have that $ |\btd \tilde g_b(u)| \le \int |\btd
f_0(y|u)| dy $ is bounded. For any $ r > 2 $, it follows that  $
\sup_{t'} E\{ \btd \tilde g_b(B_0^\top X_i) \}^r \le c$ and that
$$
\sup_{A, s, t, t'} E|m_{\ell,n}(A, s, t, t', X_i, Y_i) |^r \le c,
\quad \ell = 1, 2, 3,
$$
where $ c $ is a finite constant.  For any $ \alpha_0
>0 $, let ${{\cal D}}'_n = \{(t, t', s): |t|\le n^{\alpha_0},
|t'|\le n^{\alpha_0}, |s|\le n^{\alpha_0}\} $. By taking $ \chi =
(A, t, t', s) $ and $ a_n = c $,  we have from Lemma \ref{basic1}
\beginy
\sup_{A \in {{\cal A}}, (t,t',s) \in {{\cal D}}'_n } n^{-1}
\Big|\sum_{i=1}^n
           m_{\ell,n}(A, s, t, t', X_i, Y_i)\Big| = O( \delta_n),
           \quad \ell = 1,2,3 \label{vc1}
\endy
almost surely. On the other hand, $  |m_{\ell,n}(A, s, t, t', X_i,
Y_i)| $ is bounded. Thus,
\beginy
\sup_{A \in {{\cal A}}, (t,t',s)  } n^{-1} \Big|\sum_{i=1}^n
           m_{\ell,n}(A, s, t, t', X_i, Y_i)\Big| = O(1),
           \quad \ell = 1,2,3. \label{vc20}
\endy
By (C4), the Fourier transformation functions $ \varphi_{_K}(.) $
and $ \varphi_{_H}(.) $ are absolutely integrable; see Chung
(p.166, 1968). We can choose $ \alpha_0 $ such that
\beginy
\int_{|s|>n^{\alpha_0}}
     |\varphi_{_K}(s) |
    ds  = O(\delta_n^3), \quad
\int_{|t|>n^{\alpha_0}}
     |\varphi_{_H}(t)|dt  < O(\delta_n^3). \label{vc2}
\endy
Partition the integration region in $  \Delta_{n,1}(A) $ into two
parts, we have from (\ref{vc1})-(\ref{vc2}) that
\beginy
\n \sup_{A \in {{\cal A}}}\Big| \Delta_{n,1}(A) \Big| &\le&
\frac{1}{n^3h^{\kappa}b^2} \int_{ (s,t,t') \in {{\cal D}}_n'}
\prod_{\ell=1}^3 \sup_{A \in {{\cal A}}} | \sum_{i=1}^n
           m_{\ell,n}(A, s, t, t', X_i, Y_i)|\\
\n           && \hspace{5cm} \times |\varphi_{_K}(s)
\varphi_{_H}(t) \varphi_{_H}(t')| ds dt
              dt' \\
\n && + \frac{1}{n^3h^{\kappa}b^2} \int_{ (s,t,t') \notin {{\cal
D}}_n'} \prod_{\ell=1}^3 \sup_{A \in {{\cal A}}} | \sum_{i=1}^n
           m_{\ell,n}(A, s, t, t', X_i, Y_i)|\\
 \n          && \hspace{5cm} \times |\varphi_{_K}(s) \varphi_{_H}(t) \varphi_{_H}(t')| ds dt
              dt' \\
\n &=& (h^{\kappa}b^2)^{-1} O(\delta_n^3)\int |\varphi_{_K}(s)
\varphi_{_H}(t) \varphi_{_H}(t')| ds dt
              dt' \\
\n && \hspace{1cm} + (h^{\kappa}b^2)^{-1} O(1) \int_{ (s,t,t')
\notin {{\cal D}}_n'} |\varphi_{_K}(s) \varphi_{_H}(t)
\varphi_{_H}(t')| ds dt
              dt' \\
 &=& O( \delta_n^3 h^{-\kappa} b^{-2})  \ \  \label{ccbbdd2}
\endy
 almost surely. Let $ \tilde g(X_i) =    E_j [K_h(
A^\top X_{ij})\btd \tilde g_b(B_0^\top X_j)] $. It is easy to see
that $ \tilde g(X_i) = O(1) $ almost surely. Applying the inverse
Fourier transformation to $\varepsilon_{b,i}(Y_k) $ and using
similar arguments leading to (\ref{ccbbdd2}), we have
\beginy
\sup_{A \in {{\cal A}}} \Big| \Delta_{n,2}(A) \Big|=
 O( \delta_n^2 b^{-2}) \label{ccbbdd3}
\endy
almost surely.  Applying the inverse Fourier transformation to $
K_h( A^\top X_{ij})$, similar to (\ref{ccbbdd2}) we have
\beginy
\sup_{A \in {{\cal A}}}\Big| \Delta_{n,3}(A) \Big|  = O(
\delta_n^2 h^{-\kappa} b^{-1})  \label{ccbbdd4}
\endy
almost surely. By (\ref{cccai}), we have
\beginn
\Delta_{n,4}(A) = n^{-1}\sum_{i=1}^n E_j E_k \{ K_h( A^\top
X_{ij})  \btd  g_b(B_0^\top X_j, Y_k)\varepsilon_{b,i}(Y_k)  \}.
\endn
By Lemma \ref{basic1}, we have
\beginy
 \sup_{A \in {{\cal A}}} \Delta_{n,4}(A) = O( \delta_n) \label{ccbbdd6}
\endy
almost surely.  Finally,  Lemma \ref{basic2} follows from
(\ref{ccbbdd2})-(\ref{ccbbdd6}) and (\ref{adf07}). \hspace{\fill}
$\Box$

\

\

\noindent {\bf Acknowledgements:} Two referees and an associate
editor,  Professor Z. D. Bai and Professor B. Brown provided for
very valuable comments and suggestions for the paper. The work was
supported by NUS FRG R-155-000-048-112.\bigskip

\section*{References}

\def\ditem{\vspace{-0.2cm} \item}

\begin{description}

\ditem Ando, T. (1987) Totally positive matrices. \textit{Linear
Algebra Appl.} \textbf{90}, 165-219.

\ditem Bai, Z. D., Miao, B. Q. and Rao, C. R. (1991) Estimation of
directions of arrival of signals: Asymptotic results. {\it
Advances in Spectrum Analysis and Array Processing}, vol. 1
(edited by Simon Haykin).

\ditem  Chen, C-H. and Li, K. C. (1998) Can  SIR be as popular as
multiple linear regression? \textit{Statistica Sinica},
\textbf{8}, 289-316.


\ditem Chung K. L. (1968) \textit{A Course in probability Theory}.
Academic Press, New York.

\ditem Chow, Y. S. and Teicher, H. (1978) \textit{Probability
Theory, Independence, Interchangeability and Martingales}.
Springer-Verlag, New York.

\ditem Cook, R. D. (1998) {\it Regression Graphics}. New York,
Wiley.

\ditem Cook, R. D.  and Li, B. (2002) Dimension reduction for
conditional mean in regression.  \textit{The Annals of
Statistics}, \textbf{30}, 455-474.

\ditem Cook, R. D.  and Weisberg, S. (1991) Sliced inverse
regression for dimension reduction: comment.  \textit{J. Amer.
Statist. Ass.}, \textbf{86}, 328-332.

\ditem Delecroix, M., Hristache, M. and  Patilea, V. (2005) On
semiparametric M-estimation in single-index regression.
\textit{Journal of Statistical Planning and Inference} \textbf{
136}, 730-769

\ditem  Fan, J. and Gijbels, I. (1996) {\it Local Polynomial
Modelling and Its Applications}. Chapman $\&$ Hall, London.

\ditem Fan, J. and Yao, Q. (2003) \textit{Nonlinear Time Series :
nonparametric and parametric methods}. New York : Springer Verlag.

\ditem Fan, J., Yao, Q. and Tong, H. (1996) Estimation of
conditional densities and sensitivity measures in nonlinear
dynamical systems. {\it Biometrika}, {\bf 83}, 189-196.


\ditem H\"ardle, W. (1991) Sliced inverse regression for dimension
Reduction: comment.  \textit{J. Amer. Statist. Ass.}, \textbf{86},
332-334.

\ditem  H\"{a}rdle, W., Hall, P. and Ichimura, H. (1993) Optimal
smoothing in single-index models. {\it Ann. Statist.}, {\bf 21},
157-178.

\ditem  H\"{a}rdle, W. and Stoker, T. M. (1989) Investigating
smooth multiple regression by method of average derivatives.
\textit{J. Amer. Stat. Ass.} \textbf{84}, 986-995.

\ditem Hoeffding, W. (1961) {\it The strong law of large numbers for U-statistics}. Inst. Statist. Univ. of
       North Carolina, Mimeo Report, No. 302.

\ditem Horowitz, J. L. and Hardle, W. (1996) Direct semiparametric
estimation of single-index models with discrete covariates.
\textit{J. Amer. Stat. Ass.} \textbf{91}, 1632-1640.

\ditem Hristache, M., Juditski, A. and Spokoiny, V. (2001) Direct
estimation of the index coefficients in a single-index model. {\it
Annals of Statistics}, {\bf 29},  595-623.

\ditem Hristache, M., Juditski, A, Polzehl, J., Spokoiny, V.
(2001) Structure adaptive approach for dimension reduction. {\it
Annals of Statistics}, {\bf  29} 1537--1566.

\ditem Li, B., Zha, H. and  Chiaromonte, F. (2005) Contour
regression: a general approach to dimension reduction. {\it The
Annals of Statistics}, \textbf{33}, 1580-1616.

\ditem  Li, K. C. (1991) Sliced inverse regression for dimension
reduction (with discussion). \textit{J. Amer. Statist. Ass.},
\textbf{86}, 316-342.

\ditem  Li, K. C. (1992) On principal hessian directions for data
visualization and dimension reduction: another application of
Stein's lemma. {\it  Journal of the American Statistical
Association}, {\bf 87}, 1025-1039.

\ditem Lue, H-H. (2004) Principal Hessian directions for
regression with measurement error. \textit{Biometrika}, \textbf{
91}, 409-423.

\ditem Mack, Y. P. and  Silverman, B.W. (1982) Weak and strong
uniform consistency of kernel regression estimates. {\it Z.
Wahrsch. verw. Gebiete,} {\bf 61}, 405--415.


\ditem Samarov, A. M. (1993)  Exploring regression structure using
nonparametric functional estimation.    \textit{J. Amer. Statist.
Ass.} \textbf{88}, 836-847.

\ditem Scott, D. W. (1992) \textit{Multivariate Density
Estimation: Theory, Practice and visualization.} John Wiely \&
Sons, New York.


\ditem DE LA Pe\~na V., H. (1999) A general class of exponential
inequalities for martingales and ratios. \textit{Ann. of Prob.}
\textbf{27}, 537-564.

\ditem World Health Organization,   Reports on a WHO/HEI working
group, Bonn, Germany, (2003).

\ditem Xia, Y., Tong, H., Li, W. K. and Zhu, L. (2002) An adaptive
estimation of dimension reduction space (with discussions). {\em
J. Roy. Statist. Soc. B.}, {\bf 64}, 363-410.

\ditem Xia, Y., Tong, H. and  Li W. K. (2002) Single-index
volatility models and estimation. \textit{Statistica Sinica}, {\bf
12}, 785-799.

\ditem Xia, Y. (2006a) Asymptotic distributions for two estimators
of the single-index model. \textit{Econometric Theory}, {\bf 22},
1112-1137.

\ditem Xia, Y. (2006b)  A constructive approach to the estimation
of  dimension reduction Directions. \textit{Technical Report},
Department of Statistics and Applied Probability, National
University of Singapore.

\ditem Yin, X. and Cook, R. D. (2002) Dimension reduction for the
conditional $k$-th moment in regression. {\it J.  Roy. Stat Soc.
B}, {\bf 64}, 159-175.

\ditem Yin, X. and Cook, R. D. (2005) Direction estimation in
single-index regressions. \textit{Biometrika}, \textbf{92},
371-384.

\end{description}

\

\

 Yingcun  Xia

 Department of Statistics and Applied Probability

 National University of Singapore

 Singapore

 Email: \textit{staxyc@stat.nus.edu.sg}

\baselineskip3.0em

\end{document}